\documentclass[12pt,a4paper,final]{amsart}
\usepackage{amsrefs}
\usepackage{comment}

\numberwithin{equation}{section}

\theoremstyle{definition}
\newtheorem{definition}{Definition}[section]
\theoremstyle{plain}
\newtheorem{theorem}[definition]{Theorem}
\newtheorem{corollary}[definition]{Corollary}
\newtheorem{lemma}[definition]{Lemma}
\newtheorem{proposition}[definition]{Proposition}
\theoremstyle{remark}

\newtheorem{example}[definition]{Example}

\def\beq{\begin{equation}}
\def\eeq{\end{equation}}
\def\bproof{\begin{proof}}
\def\eproof{\end{proof}}
\def\pl{\partial}
\def\R{\mathbb{R}}

\def\ol{\overline}
\def\tim{\times}  
\def\gb{\beta}  
\def\mid{\,:\,}   
\def\bcases{\begin{cases}}
\def\ecases{\end{cases}}
\def\gep{\varepsilon}      
\def\ep{\gep} 
\def\gl{\lambda}
\def\gam{\gamma}
\def\fr{\frac}
\def\gth{\theta}
\def\gd{\delta}
\def\go{\omega}
\def\bye{\end{document}}

\def\bmat{\begin{pmatrix}}
\def\emat{\end{pmatrix}}

\def\bS{\mathbb{S}}
\DeclareMathOperator{\cone}{cone}
\DeclareMathOperator*{\einf}{ess\, inf}
\DeclareMathOperator*{\esup}{ess\, sup}
\DeclareMathOperator*{\elim}{ess\, lim}

\DeclareMathOperator{\diag}{diag}

\def\IN{\text{ in } }
\def\AND{\text{ and }}
\def\FOR{\text{ for }}
\def\FORALL{\text{ for all }}
\def\ON{\text{ on }}
\def\IF{\text{ if }}
\def\bald{\begin{aligned}}
\def\eald{\end{aligned}}
\def\stm{\setminus}
\def\erf{\eqref}

\def\ga{\alpha}
\def\gs{\sigma}

\def\gO{\Omega}

\def\AE{\text{ a.e. }}

\def\tH{\widetilde H}

\title{Nonlinear Neumann problems for fully nonlinear elliptic PDEs on a quadrant}

\author[H. Ishii]{Hitoshi Ishii$^*$}
\address[H. Ishii]{Institute for Mathematics and Computer Science, Tsuda  University,
 2-1-1 Tsuda, Kodaira, Tokyo 187-8577 Japan.}
\email{hitoshi.ishii@waseda.jp}
\thanks{${}^*$ Corresponding author}

\author[T. Kumagai]{Taiga Kumagai}
\address[T. Kumagai]{National Institute of Technology, Maizuru College, 234 Shiroya, Maizuru-shi, Kyoto 625-8511 Japan.}
\email{kumatai13@gmail.com}

\date{\today}


\thanks{H. Ishii was partially supported by the JSPS KAKENHI Grant Nos. JP16H03948, JP20K03688, JP20H01817 and JP21H00717; T. Kumagai 
was partially supported by the JSPS KAKENHI Grant No. JP20K03688.}


\date{\today}

\begin{document}

\begin{abstract} We consider the nonlinear Neumann problem for fully nonlinear 
elliptic PDEs on a quadrant. We establish a comparison theorem for viscosity sub and supersolutions of the nonlinear Neumann problem. The crucial argument in the proof 
of the comparison theorem is to build a $C^{1,1}$ test function which takes care of 
the nonlinear Neumann boundary condition. 
A similar problem has been treated on a general $n$-dimensional orthant by Biswas, Ishii, Subhamay, and Wang 
[SIAM J. Control Optim. 55 (2017), pp. 365--396], where the functions 
($H_i$ in the main text) describing 
the boundary condition are required to be positively one-homogeneous, 
and the result in this paper removes the positive homogeneity in two-dimension. 
An existence result for solutions is also presented. 
\end{abstract}

\maketitle

\tableofcontents

\allowdisplaybreaks

\section{Introduction}
We consider the nonlinear Neumann boundary value problem 
\beq\label{N-prob}\left\{\bald
&F(D^2u,Du,u,x)=0 \ \ \IN \gO,
\\&H_1(Du)=0 \ \ \ON \pl_1\gO,\quad H_2(Du)=0 \ \ \ON \pl_2\gO,
\eald \right. 
\eeq
where $\gO=\{x=(x_1,x_2)\in\R^2\mid x_1>0,\, x_2>0 \}$ (called the first quadrant), 
$\pl_1\gO=\{x=(x_1,x_2)\in\R^2\mid x_1=0,\,x_2\geq 0\}$, and 
$\pl_2\gO=\{x=(x_1,x_2)\in\R^2\mid x_1\geq 0,\,x_2=0\}$. Furthermore, $F$ is a fully
nonlinear, degenerate elliptic operator for functions $u$ on $\ol\gO$, and the Neumann type operators $H_i$ are given as follows. 
\beq\tag{B1}\label{cond-H_i} \ \left\{\ \bald
& \text{For } i=1,2, \ H_i(p)=-p_i+h_i(p_j), \ \ \text{ where } j=2 \IF i=1 \AND j=1 \IF i=2,
\\& h_i \ \text{ convex}, 
\\ & h_i'(t)\geq 0\ \ \AE t\in\R, \quad \elim_{t\to \infty}h_i'(t)=\ga_i\in[0,\,1)
. 
\eald \right.
\eeq

We remark that, since $h_i$ above are convex on $\R$, $h_i$ are differentiable 
almost everywhere and $h_i'$ are ``nondecreasing" on $\R$, and  
we have \ $\elim_{t\to \infty}h'_i(t)=\esup_{t\in \R}h_i'(t)$. Similarly, we have \ 
$\elim_{t\to 0+} h_i'(t)=\einf_{t>0}h_i'(t)$ and  
$\elim_{t\to 0-} h_i'(t)=\esup_{t<0}h_i'(t)$, for instance. 

The condition that 
$h_i'(t)\geq 0$ a.e. on $\R$ can be stated that $h_i$ is nondecreasing on 
$\R$. 

In \cite{BISW2017}, the authors have studied the nonlinear elliptic partial 
differential equation
\beq\label{n-dPDE}
F(D^2u,Du,u,x)=0 \ \ \IN  \gO:=(0,\infty)^n,
\eeq 
with the nonlinear Neumann type 
boundary condition, similar to the above, stated as
\beq\label{n-dBC}
H_i(Du)=0 \ \ \ON \pl_i \gO:=\{x=(x_1,\ldots,x_n)\in [0,\infty)^n\mid x_i=0\},
\eeq
where \ $H_i(p):=-p_i+h_i(p).$ 

Henceforth, let $\{e_1,\ldots,e_n\}$ denote the standard basis of $\R^n$
and the symbol $\,\cdot\,$ indicate the Euclidean inner product in $\R^n$. 
For each $i\in\mathbb{I}:=\{1,\ldots,n\}$, the function $h_i$, defining \erf{n-dBC}, 
has the form
\beq\label{n-dh}
h_i(p):=\max_{R\in \mathbf{R}_i} (p\cdot Re_i-g_i(R))=\max_{R\in \mathbf{R}_i} 
\left(\sum_{j\in \mathbb{I}\stm\{i\}}p_jR_{ji}-g_i(R)\right),
\eeq
where $\mathbf{R}_i\subset \mathbf{R}$, with $\mathbf{R}$ denoting   
the set of $n\tim n$ real matrices $R=(R_{jk})\in \R^{n\tim n}$ such that  $R_{jj}=0$ and $R_{jk}\geq 0$ for all $j,k\in\mathbb{I}$, 
\[
\mathbf{R}_i:=\{R\in\mathbf{R}\mid \sum_{j\in \mathbb{I}}R_{ji}\leq \gb_i\},
\] 
for a given $\gb_i\in [0,1)$, and $g_i\in C(\mathbf{R}_i)$ is a given function. 
Observe that the functions $h_i$ are convex functions on $\R^n$ and that for any $i\in\mathbb{I}$,
\[
p\cdot Re_i=\sum_{j\in \mathbb{I}\stm\{i\}}p_jR_{ji}
\]
and hence, $h_i(p)$ is independent of $p_i$. For any $p,q\in\R^n$, we have 
\begin{align*}
&|p\cdot Re_i-q\cdot Re_i|\leq 
|p-q||Re_i|,
\\&|Re_i|^2 =\sum_{j\in\mathbb{I}}R_{ji}^2\leq \gb_i\sum_{j\in\mathbb{I}}R_{ji}\leq \gb_i^2, \\
\intertext{and hence, }  
&
|p\cdot Re_i-q\cdot Re_i|\leq \gb_i |p-q|,
\end{align*}
which implies that $h_i$ is Lipschitz continuous on $\R^n$ with $\gb_i$ as a Lipschitz 
bound. 
When $n=2$, with a slight abuse of notation, we may write $h_1(p)=h_1(p_2)$, $h_2(p)=h_2(p_1)$. The functions $h_i$ are obviously nondecreasing on $\R$ by \erf{n-dh}, and satisfy $\esup_{\R} h_i'\in [0,\gb_i]$. Thus, condition \erf{cond-H_i} 
is satisfied in this case. 

As shown in \cite{BISW2017}, the boundary value problem \erf{n-dPDE}--\erf{n-dBC} arises in a scaling limit of stochastic control of queuing systems. 
A main contribution in \cite{BISW2017} is that, under suitable hypotheses, 
the value function of the stochastic control problem, obtained in the scaling limit, 
is identified as the unique viscosity solution of \erf{n-dPDE} and \erf{n-dBC}, 
with an appropriate choice of $F$. This identification result heavily depends 
on the uniqueness or comparison theorem for viscosity solutions, 
subsolutions, and supersolutions of \erf{n-dPDE}--\erf{n-dBC}. 
The comparison theorem in \cite{BISW2017} is basically valid only for the $H_i$ having 
the positive homogeneity of degree one.  

The aim of this paper is to attempt to remove the above requirement that the $H_i$ should be positively homogeneous of degree one. We are successful to remove the homogeneity requirement on the $H_i$, but only in dimension two. The comparison principle 
for \erf{n-dPDE}--\erf{n-dBC}, with $h_i$ given by \erf{n-dh}, 
is still largely open in higher dimensions and even in dimension two (see \erf{cond-slope}).   

We remark that there are a great amount of contributions  to the Neumann 
boundary value problem on relatively regular domains even in the framework of 
viscosity solutions. We refer for some general results to
\cites{CIL1992, I1991, B1999}, and for some results on domains with corners
to \cites{DI1991, BISW2017} and the references therein.

The crucial part of the proof of our comparison principle is to build an appropriate 
convex test function $f$ on $\R^2$.  As described in Corollary \ref{exist-f} below, 
the required properties on $f$ are that $f$ is a $C^{1,1}$ function, has a 
quadratic growth, and  satisfies 
\beq\label{cond-f}
H_1(Df(x))
\bcases 
\geq 0\ \ &\IF\ x_1\leq 0, \\
\leq 0 &\IF\ x_1\geq 0, 
\ecases
\quad\AND \quad
H_2(Df(x))
\bcases
\geq 0 \ \ &\IF\ x_2\leq 0, \\
\leq 0 &\IF\ x_2\geq 0.
\ecases
\eeq
Here and henceforth, we say that a function $g$ defined on an open subset $U$ of $\R^2$ is of class $C^{1,1}$ if the first derivatives of $g$ is bounded on $U$. 

Our strategy of building such a test function $f$, a version of which 
appears as $\phi_\ep$ 
in Section 5,  is first to define 
a convex function $\psi$ (see \erf{def-psi} below) on $\R^2$ associated with the curves $H_i=0$, $i=1,2$,
having a linear growth, and then to take the convex conjugate of the function 
$\psi^2$ plus a small quadratic function $q_\ep\mid x\mapsto \ep|x|^2$, with 
$\ep>0$.  Roughly speaking, the addition of $q_\ep$ before taking the convex conjugate corresponds to the regularization by the inf-convolution 
(called also the Yosida approximation) of the convex function $\psi^2$.  
See \erf{inf-conv} below for this regularization.)

The paper is organized as follows. Let $f$, $\psi$, and $q_\ep$ be the functions 
appearing in the explanation above. Throughout Sections 2--4, we 
assume that $H_i$, $i=1,2$, vanish at the origin. This normalization condition is 
called \erf{zero-norm} below, and is removed in Section 6. Under the normalization assumption, $f$ is called $\phi_\ep$ (see \erf{def-phi} below). Section 2 is devoted to the study of 
the geometry of the curves $H_i=0$, $i=1,2$, and certain line segments 
connecting points on those curves, which will be the level sets of the function $\psi$. 
Section 3 introduces a condition (see \erf{cond-slope} below), which is 
crucial to define the function $\psi$, with necessary properties, together with three 
examples of cases of $H_i$, $i=1,2$, for which the critical condition \erf{cond-slope} is valid. Section 4 is devoted to building 
the function $\psi$ and the study of its basic properties.  
We define the function $\phi_\ep$ as the inf-convolution of the convex 
$\psi^2+q_\ep$ in Section 5. In Section 6, the function $f$ is defined as the perturbation of
 $\phi_\ep$ by a linear function, which fits the general case, and we establish a comparison theorem for sub ans supersolutions of \erf{N-prob} as an application of the test function $f$. An existence result for bounded solutions of \erf{N-prob} is presented as well.

\section{Some geometric preliminaries} 
We introduce a normalization assumption. 
\beq\tag{B2}\label{zero-norm} 
h_i(0)=0 \ \ \FOR i=1,2.
\eeq
This condition is clearly equivalent to the condition that $H_i(0)=0$ for $i=1,2$.  
We \emph{assume} throughout this section that \erf{cond-H_i}--\erf{zero-norm} are satisfied. We come back to the general case without \erf{zero-norm} in Section 6. 
We set 
\begin{align*}
P&=\{x=(x_1,x_2)\in\R^2\mid x_1\geq h_1(x_2),\, x_2\geq h_2(x_1)\},
\\ Q&=\{x\in\R^2\mid x_1\leq h_1(x_2),\, x_2\geq h_2(x_1)\}, 
\\ R&=\{x\in\R^2\mid x_1\geq h_1(x_2),\, x_2\leq h_2(x_1)\}, 
\\ S&=\{x\in\R^2\mid x_1\leq h_1(x_2),\, x_2\leq h_2(x_1)\}.
\end{align*}
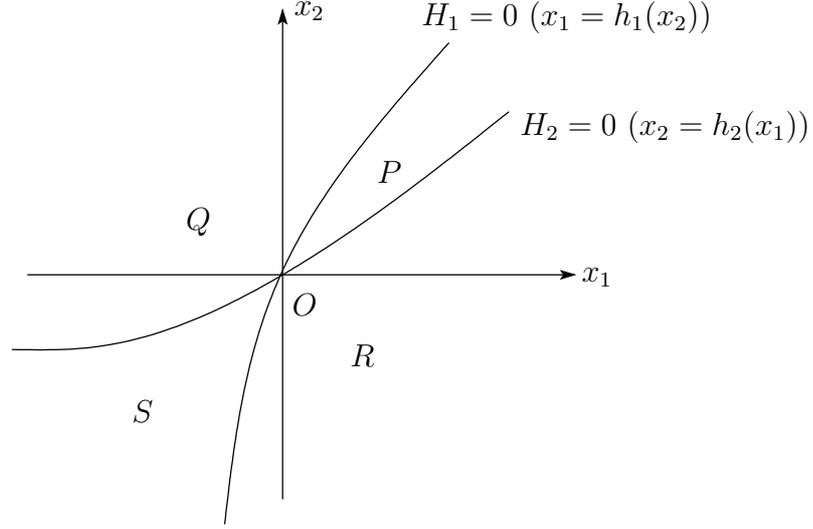
\begin{center}
{\unitlength 0.1in%
\begin{picture}(29.5000,27.6000)(8.2000,-29.0000)%
%
\special{pn 8}%
\special{pa 1920 2900}%
\special{pa 1924 2868}%
\special{pa 1927 2836}%
\special{pa 1939 2740}%
\special{pa 1942 2709}%
\special{pa 1954 2613}%
\special{pa 1959 2581}%
\special{pa 1963 2549}%
\special{pa 1967 2518}%
\special{pa 1977 2454}%
\special{pa 1981 2423}%
\special{pa 1987 2391}%
\special{pa 1992 2359}%
\special{pa 1997 2328}%
\special{pa 2003 2296}%
\special{pa 2015 2234}%
\special{pa 2021 2202}%
\special{pa 2042 2109}%
\special{pa 2049 2077}%
\special{pa 2065 2015}%
\special{pa 2073 1985}%
\special{pa 2100 1892}%
\special{pa 2110 1862}%
\special{pa 2120 1831}%
\special{pa 2130 1801}%
\special{pa 2152 1741}%
\special{pa 2164 1711}%
\special{pa 2175 1681}%
\special{pa 2188 1652}%
\special{pa 2200 1622}%
\special{pa 2226 1564}%
\special{pa 2240 1535}%
\special{pa 2254 1507}%
\special{pa 2268 1478}%
\special{pa 2313 1394}%
\special{pa 2329 1366}%
\special{pa 2345 1339}%
\special{pa 2361 1311}%
\special{pa 2377 1284}%
\special{pa 2411 1230}%
\special{pa 2429 1204}%
\special{pa 2446 1177}%
\special{pa 2464 1150}%
\special{pa 2500 1098}%
\special{pa 2519 1072}%
\special{pa 2537 1046}%
\special{pa 2575 994}%
\special{pa 2595 969}%
\special{pa 2614 943}%
\special{pa 2634 918}%
\special{pa 2654 892}%
\special{pa 2734 792}%
\special{pa 2776 742}%
\special{pa 2796 717}%
\special{pa 2817 692}%
\special{pa 2838 668}%
\special{pa 2880 618}%
\special{pa 2902 594}%
\special{pa 2923 569}%
\special{pa 2944 545}%
\special{pa 2966 520}%
\special{pa 2987 496}%
\special{pa 3008 471}%
\special{pa 3030 447}%
\special{pa 3051 422}%
\special{pa 3073 398}%
\special{pa 3080 390}%
\special{fp}%
%
\special{pn 8}%
\special{pa 820 1990}%
\special{pa 852 1991}%
\special{pa 916 1991}%
\special{pa 949 1992}%
\special{pa 981 1992}%
\special{pa 1013 1991}%
\special{pa 1045 1991}%
\special{pa 1109 1989}%
\special{pa 1141 1987}%
\special{pa 1172 1984}%
\special{pa 1204 1982}%
\special{pa 1236 1978}%
\special{pa 1267 1974}%
\special{pa 1299 1969}%
\special{pa 1330 1964}%
\special{pa 1392 1952}%
\special{pa 1424 1945}%
\special{pa 1455 1937}%
\special{pa 1485 1929}%
\special{pa 1516 1921}%
\special{pa 1547 1912}%
\special{pa 1578 1902}%
\special{pa 1608 1892}%
\special{pa 1638 1881}%
\special{pa 1669 1871}%
\special{pa 1789 1823}%
\special{pa 1818 1810}%
\special{pa 1848 1797}%
\special{pa 1878 1783}%
\special{pa 1936 1755}%
\special{pa 1965 1740}%
\special{pa 1994 1726}%
\special{pa 2023 1710}%
\special{pa 2081 1680}%
\special{pa 2137 1648}%
\special{pa 2166 1632}%
\special{pa 2194 1615}%
\special{pa 2222 1599}%
\special{pa 2250 1582}%
\special{pa 2277 1565}%
\special{pa 2305 1549}%
\special{pa 2332 1531}%
\special{pa 2360 1514}%
\special{pa 2387 1497}%
\special{pa 2414 1479}%
\special{pa 2441 1462}%
\special{pa 2495 1426}%
\special{pa 2521 1408}%
\special{pa 2548 1390}%
\special{pa 2574 1372}%
\special{pa 2601 1354}%
\special{pa 2627 1335}%
\special{pa 2653 1317}%
\special{pa 2679 1298}%
\special{pa 2705 1280}%
\special{pa 2809 1204}%
\special{pa 2834 1185}%
\special{pa 2860 1166}%
\special{pa 2885 1147}%
\special{pa 2911 1127}%
\special{pa 2936 1108}%
\special{pa 2962 1089}%
\special{pa 2987 1069}%
\special{pa 3012 1050}%
\special{pa 3062 1010}%
\special{pa 3088 991}%
\special{pa 3163 931}%
\special{pa 3188 912}%
\special{pa 3238 872}%
\special{pa 3262 852}%
\special{pa 3362 772}%
\special{pa 3387 753}%
\special{pa 3390 750}%
\special{fp}%
\put(34.5000,-9.1000){\makebox(0,0)[lb]{{$H_2=0$ $(x_2=h_2(x_1))$}}}%
\put(29.4000,-3.4000){\makebox(0,0)[lb]{{$H_1=0$ $(x_1=h_1(x_2))$}}}%
\put(22.7000,-18.1000){\makebox(0,0)[lb]{{$O$}}}%
\put(27.1000,-11.2000){\makebox(0,0)[lb]{{$P$}}}%
\put(17.2000,-14.0000){\makebox(0,0)[lb]{{$Q$}}}%
\put(25.7000,-20.8000){\makebox(0,0)[lb]{$R$}}%
\put(14.4000,-23.7000){\makebox(0,0)[lb]{$S$}}%
\put(37.7000,-16.6000){\makebox(0,0)[lb]{$x_1$}}%
\put(22.8000,-2.7000){\makebox(0,0)[lb]{{$x_2$}}}%
%
\special{pn 8}%
\special{pa 2220 2770}%
\special{pa 2220 220}%
\special{fp}%
\special{sh 1}%
\special{pa 2220 220}%
\special{pa 2200 287}%
\special{pa 2220 273}%
\special{pa 2240 287}%
\special{pa 2220 220}%
\special{fp}%
%
\special{pn 8}%
\special{pa 900 1600}%
\special{pa 3730 1600}%
\special{fp}%
\special{sh 1}%
\special{pa 3730 1600}%
\special{pa 3663 1580}%
\special{pa 3677 1600}%
\special{pa 3663 1620}%
\special{pa 3730 1600}%
\special{fp}%
\end{picture}}%

\end{center}

We fix $v_0=(v_{01},v_{02}),v_1=(v_{11}.v_{12}),v_2=(v_{21},v_{22})\in\R^2$ such that  
\begin{align}
\label{cond-v_i}& v_0\cdot e_1>0,\ v_0\cdot e_2>0,\quad 
v_1\cdot (-e_1)>0,\ v_1\cdot e_2>0,\quad v_2\cdot e_1>0, \quad v_2\cdot (-e_2)>0,
\\\label{cond-v_1}&\einf_{t\geq 0}v_1\cdot(h_1'(t),1)>0,\quad \einf_{t\leq 0}v_1\cdot(-1,-h_2'(t))>0,
\\\label{cond-v_2}&\einf_{t\geq 0}v_2\cdot(1, h_2'(t))>0,\quad \einf_{t\leq 0}v_2\cdot (-h_1'(t),-1)>0.
\end{align}
Note that, since $h_1'(t)\geq 0$ and $h_2'(t)\geq 0$ a.e.,  
\beq \label{cond-v_0}\einf_{t\geq 0}v_0\cdot (1,h_2'(t))\geq v_{01}>0,\quad \einf_{t\geq 0}v_0\cdot (h_1'(t),1)\geq v_{02}>0.\eeq
Since $0\leq h_i'(t)<1$ a.e. on $\R$ for $i=1,2$, it is easily seen that 
\erf{cond-v_1} (resp., \erf{cond-v_2}) implies that 
$v_1\cdot e_1<0$ and $v_1\cdot e_2>0$ (resp., $v_2\cdot e_1>0$ and $v_2\cdot e_2<0$).  Moreover, we deduce easily that \erf{cond-v_1} (resp., \erf{cond-v_2}) 
is equivalent to the condition that $v_1\cdot(\ga_1,1)>0$ and 
$v_1\cdot (-1,-\ga_2^0)>0$ (resp., $v_2\cdot (1,\ga_2)>0$ and $v_2\cdot (-\ga_1^0,-1)>0$), where $\ga_i^0=\esup_{t\leq 0}h_i'(t)$ for $i=1,2$.

Let $\gl>0$. Let $x_{PQ}^\gl,\, x_{PR}^\gl$ be the points 
where the line $v_0\cdot x=\gl$ intersects the curves $H_1=0$ and 
$H_2=0$, respectively. Set $\gl_1=v_1\cdot x_{PQ}^\gl$ and 
$\gl_2=v_2\cdot x_{PR}^\gl$.  Note that $\gl_1$ and $\gl_2$ are positive constants.  
Let $x_{QS}^\gl,\, x_{RS}^\gl$ be the points in $(-\infty,0)^2$
where the lines $v_1\cdot x=\gl_1$ and $v_2\cdot x=\gl_2$ intersect, 
respectively, the curves $H_2=0$ and $H_1=0$. The existence and uniqueness of 
$(x_{PQ}^\gl,x_{PR}^\gl,x_{QS}^\gl,x_{RS}^\gl)$ is established in the next proposition. 

\begin{center}
{\unitlength 0.1in%
\begin{picture}(31.1000,28.8000)(7.9000,-29.9000)%
%
\special{pn 8}%
\special{pa 790 1600}%
\special{pa 3840 1600}%
\special{fp}%
\special{sh 1}%
\special{pa 3840 1600}%
\special{pa 3773 1580}%
\special{pa 3787 1600}%
\special{pa 3773 1620}%
\special{pa 3840 1600}%
\special{fp}%
\special{pa 2210 2990}%
\special{pa 2210 120}%
\special{fp}%
\special{sh 1}%
\special{pa 2210 120}%
\special{pa 2190 187}%
\special{pa 2210 173}%
\special{pa 2230 187}%
\special{pa 2210 120}%
\special{fp}%
%
\special{pn 8}%
\special{pa 1920 2900}%
\special{pa 1924 2868}%
\special{pa 1927 2836}%
\special{pa 1939 2740}%
\special{pa 1942 2709}%
\special{pa 1954 2613}%
\special{pa 1959 2581}%
\special{pa 1963 2549}%
\special{pa 1967 2518}%
\special{pa 1977 2454}%
\special{pa 1981 2423}%
\special{pa 1987 2391}%
\special{pa 1992 2359}%
\special{pa 1997 2328}%
\special{pa 2003 2296}%
\special{pa 2015 2234}%
\special{pa 2021 2202}%
\special{pa 2042 2109}%
\special{pa 2049 2077}%
\special{pa 2065 2015}%
\special{pa 2073 1985}%
\special{pa 2100 1892}%
\special{pa 2110 1862}%
\special{pa 2120 1831}%
\special{pa 2130 1801}%
\special{pa 2152 1741}%
\special{pa 2164 1711}%
\special{pa 2175 1681}%
\special{pa 2188 1652}%
\special{pa 2200 1622}%
\special{pa 2226 1564}%
\special{pa 2240 1535}%
\special{pa 2254 1507}%
\special{pa 2268 1478}%
\special{pa 2313 1394}%
\special{pa 2329 1366}%
\special{pa 2345 1339}%
\special{pa 2361 1311}%
\special{pa 2377 1284}%
\special{pa 2411 1230}%
\special{pa 2429 1204}%
\special{pa 2446 1177}%
\special{pa 2464 1150}%
\special{pa 2500 1098}%
\special{pa 2519 1072}%
\special{pa 2537 1046}%
\special{pa 2575 994}%
\special{pa 2595 969}%
\special{pa 2614 943}%
\special{pa 2634 918}%
\special{pa 2654 892}%
\special{pa 2734 792}%
\special{pa 2776 742}%
\special{pa 2796 717}%
\special{pa 2817 692}%
\special{pa 2838 668}%
\special{pa 2880 618}%
\special{pa 2902 594}%
\special{pa 2923 569}%
\special{pa 2944 545}%
\special{pa 2966 520}%
\special{pa 2987 496}%
\special{pa 3008 471}%
\special{pa 3030 447}%
\special{pa 3051 422}%
\special{pa 3073 398}%
\special{pa 3080 390}%
\special{fp}%
%
\special{pn 8}%
\special{pa 820 1990}%
\special{pa 852 1991}%
\special{pa 916 1991}%
\special{pa 949 1992}%
\special{pa 981 1992}%
\special{pa 1013 1991}%
\special{pa 1045 1991}%
\special{pa 1109 1989}%
\special{pa 1141 1987}%
\special{pa 1172 1984}%
\special{pa 1204 1982}%
\special{pa 1236 1978}%
\special{pa 1267 1974}%
\special{pa 1299 1969}%
\special{pa 1330 1964}%
\special{pa 1392 1952}%
\special{pa 1424 1945}%
\special{pa 1455 1937}%
\special{pa 1485 1929}%
\special{pa 1516 1921}%
\special{pa 1547 1912}%
\special{pa 1578 1902}%
\special{pa 1608 1892}%
\special{pa 1638 1881}%
\special{pa 1669 1871}%
\special{pa 1789 1823}%
\special{pa 1818 1810}%
\special{pa 1848 1797}%
\special{pa 1878 1783}%
\special{pa 1936 1755}%
\special{pa 1965 1740}%
\special{pa 1994 1726}%
\special{pa 2023 1710}%
\special{pa 2081 1680}%
\special{pa 2137 1648}%
\special{pa 2166 1632}%
\special{pa 2194 1615}%
\special{pa 2222 1599}%
\special{pa 2250 1582}%
\special{pa 2277 1565}%
\special{pa 2305 1549}%
\special{pa 2332 1531}%
\special{pa 2360 1514}%
\special{pa 2387 1497}%
\special{pa 2414 1479}%
\special{pa 2441 1462}%
\special{pa 2495 1426}%
\special{pa 2521 1408}%
\special{pa 2548 1390}%
\special{pa 2574 1372}%
\special{pa 2601 1354}%
\special{pa 2627 1335}%
\special{pa 2653 1317}%
\special{pa 2679 1298}%
\special{pa 2705 1280}%
\special{pa 2809 1204}%
\special{pa 2834 1185}%
\special{pa 2860 1166}%
\special{pa 2885 1147}%
\special{pa 2911 1127}%
\special{pa 2936 1108}%
\special{pa 2962 1089}%
\special{pa 2987 1069}%
\special{pa 3012 1050}%
\special{pa 3062 1010}%
\special{pa 3088 991}%
\special{pa 3163 931}%
\special{pa 3188 912}%
\special{pa 3238 872}%
\special{pa 3262 852}%
\special{pa 3362 772}%
\special{pa 3387 753}%
\special{pa 3390 750}%
\special{fp}%
\put(34.5000,-9.1000){\makebox(0,0)[lb]{$H_2=0$}}%
\put(31.3000,-4.6000){\makebox(0,0)[lb]{$H_1=0$}}%
\put(22.7000,-18.1000){\makebox(0,0)[lb]{$O$}}%
%
\special{pn 8}%
\special{pa 2860 1020}%
\special{pa 2990 750}%
\special{fp}%
\special{sh 1}%
\special{pa 2990 750}%
\special{pa 2943 801}%
\special{pa 2967 798}%
\special{pa 2979 819}%
\special{pa 2990 750}%
\special{fp}%
%
\special{pn 8}%
\special{pa 2650 920}%
\special{pa 2960 1080}%
\special{fp}%
\special{pa 2960 1080}%
\special{pa 2960 1080}%
\special{fp}%
%
\special{pn 8}%
\special{pa 2970 1080}%
\special{pa 1960 2670}%
\special{fp}%
%
\special{pn 8}%
\special{pa 2630 930}%
\special{pa 1350 1950}%
\special{fp}%
\special{pa 1340 1960}%
\special{pa 1340 1960}%
\special{fp}%
%
\special{pn 8}%
\special{pa 2420 1940}%
\special{pa 2750 2150}%
\special{fp}%
\special{sh 1}%
\special{pa 2750 2150}%
\special{pa 2704 2097}%
\special{pa 2705 2121}%
\special{pa 2683 2131}%
\special{pa 2750 2150}%
\special{fp}%
%
\special{pn 8}%
\special{pa 1970 1460}%
\special{pa 1710 1170}%
\special{fp}%
\special{sh 1}%
\special{pa 1710 1170}%
\special{pa 1740 1233}%
\special{pa 1746 1210}%
\special{pa 1769 1206}%
\special{pa 1710 1170}%
\special{fp}%
\put(30.4000,-7.4000){\makebox(0,0)[lb]{$v_0$}}%
\put(27.8000,-21.5000){\makebox(0,0)[lb]{$v_2$}}%
\put(17.4000,-11.4000){\makebox(0,0)[lb]{$v_1$}}%
\put(29.9000,-12.2000){\makebox(0,0)[lb]{$x_{PR}^\lambda$}}%
\put(22.8000,-8.6000){\makebox(0,0)[lb]{$x_{PQ}^\lambda$}}%
\put(11.3000,-22.2000){\makebox(0,0)[lb]{$x_{QS}^\lambda$}}%
\put(15.9000,-27.5000){\makebox(0,0)[lb]{$x_{RS}^\lambda$}}%
\put(11.2000,-13.6000){\makebox(0,0)[lb]{$Q$}}%
\put(25.1000,-24.2000){\makebox(0,0)[lb]{$R$}}%
\put(14.4000,-25.0000){\makebox(0,0)[lb]{$S$}}%
\put(39.0000,-17.2000){\makebox(0,0)[lb]{$x_1$}}%
\put(22.9000,-2.4000){\makebox(0,0)[lb]{$x_2$}}%
\end{picture}}%
\end{center}

\begin{proposition} \label{four-points} Under \erf{cond-H_i}--\erf{zero-norm}, 
for each $\gl>0$, there exists a unique $(x_{PQ}^\gl,x_{PR}^\gl,x_{QS}^\gl,x_{RS}^\gl)\in (\R^2)^4$ with the properties described above. 
\end{proposition}

Let $I$ be an interval of $\R$ and $f\mid I\to\R$. 
In this paper, we say that $f$ is \emph{uniformly increasing} on $I$
if $t\mapsto f(t)-\ep t$ is nondecreasing on $I$ for some constant $\ep>0$. 
Also, we say that $f$ is \emph{uniformly decreasing} on $I$
if $t\mapsto f(t)+\ep t$ is nonincreasing on $I$ for some constant $\ep>0$. 
It is easily seen that if $f\mid I \to\R$ is uniformly increasing (resp., decreasing) and 
Lipschitz continuous on $I$, then $f\mid I\to f(I)$ has an inverse function 
$f^{-1}$ which is uniformly increasing (resp., decreasing) and Lipschitz continuous on $f(I)$.   

We introduce the functions $f_i$, $i=1,2$, on $[0,\infty)$,
$g_i$, $i=1,2$, on $[0,\infty)$, $p_i$, $i=1,2$, on $(-\infty,0]$, 
$q_i$, $i=1,2$, on $(-\infty,0]$ by setting   
\beq \label{fgpq} \left\{\ \bald
&f_1(t)=v_0\cdot(h_1(t),t), \ \ f_2(t)=v_0\cdot(t,h_2(t)), \ \ 
\\&g_1(t)=v_1\cdot (h_1(t),t),\ \ g_2(t)=v_2\cdot (t,h_2(t)), \ \
\\& p_1(t)=v_1\cdot (t,h_2(t)), \ \ p_2(t)=v_2\cdot (h_1(t),t),
\\& q_1(t)=m_\gam\cdot (t,h_2(t)), \ \ q_2(t)=n_\gam\cdot (h_1(t),t),
\eald \right. 
\eeq
where $m_\gam:=(-1,-\gam)$ and $n_\gam:=(-\gam,-1)$, with $\gam$ being 
any fixed positive constant. In the later sections, the constant $\gam$ will be fixed 
as the one from \erf{cond-slope}.

\begin{lemma} \label{4 func} The functions $f_i$ and $g_i$, with $i=1,2$, are uniformly 
increasing and Lipschitz continuous on $[0,\infty)$. They have inverse functions 
$f_i^{-1}$ and $g_i^{-1}$, defined on $[0,\infty)$, and the inverse functions 
$f_i^{-1},\, g_i^{-1}$ are uniformly increasing and Lipschitz continuous. 
Similarly, the functions $p_i, q_i$, $i=1,2$,  are uniformly 
decreasing and Lipschitz continuous on $(-\infty,0]$ and they have inverse functions 
$p_i^{-1}$ and $q_i^{-1}$, defined on $[0,\infty)$, and the inverse functions 
$p_i^{-1},\, q_i^{-1}$ are uniformly decreasing and Lipschitz continuous. 
\end{lemma} 

\bproof By \erf{cond-H_i}, the functions $h_i$ are Lipschitz continuous on $\R$ and hence, the functions $f_i,\,g_i,p_i,q_i$, $i=1,2$, are Lipschitz continuous. Note also that 
$f_i(0)=g_i(0)=p_i(0)=q_i(0)=0$ for all $i=1,2$.

We have for a.e. $t\geq 0$,
\[\bald
&f_1'(t)=v_0\cdot (h_1'(t),1),\quad f_2'(t)=v_0\cdot (1,h_2'(t)), \quad 
\\&g_1'(t)=v_1\cdot (h_1'(t),1), \quad   g_2'(t)=v_2\cdot(1,h_2'(t)). 
\eald
\] 
By \erf{cond-v_0}, \erf{cond-v_1} and \erf{cond-v_2}, we find that for some 
constant $\ep_1>0$,
\[
f'_i(t)\geq \ep_1,\ g'_i(t)\geq \ep_1 \ \ \text{ for a.e. }t\geq 0 \text{ and all }i=1,2,
\]
which show that $f_i,\,g_i$, $i=1,2$, are uniformly increasing 
and map onto $[0,\infty)$.  
Now, since the functions $f_i,\,g_i\mid [0,\infty)\to[0,\infty)$, $i=1,2$, are 
uniformly increasing, Lipschitz continuous, and surjective, we conclude as well 
that their inverse functions exist and are uniformly increasing and Lipschitz continuous.

Next, note that for a.e. $t\leq 0$,
\[
m_\gam\cdot (1,h_2'(t))=-1-\gam h_2'(t)\leq -1 \ \ \AND \ \ 
n_\gam\cdot (h_1'(t),1)=-\gam h_1'(t)-1\leq -1.  
\]
Using these,  \erf{cond-v_1}, and \erf{cond-v_2}, we easily find that for a.e. $t\leq 0$, all $i=1,2$, and some constant $\ep_2>0$,
\[
p_i'(t)\leq -\ep_2 \ \ \AND \ \ q_i'(t)\leq -1,
\]
which show that the functions $p_i, q_i$, $i=1,2$, are uniformly decreasing 
and map onto $[0,\infty)$. Arguing as above, we conclude that the functions 
$p_i, q_i\mid(-\infty,0]\to[0,\infty)$ are uniformly decreasing and Lipschitz continuous and have their inverse $p_i^{-1}, q_i^{-1}$ and that $p_i^{-1}, q_i^{-1}$ are 
uniformly decreasing and Lipschitz continuous on $[0,\infty)$.  The proof is complete. 
\eproof

In what follows, we use the functions $G_i\mid \R\to \R^2$, $i=1,2$, defined by  
\[
G_1(t)=(h_1(t),t) \ \ \AND  \ \ G_2(t)=(t,h_2(t)).
\]

\bproof[Proof of Proposition \ref{four-points}]  

{\sc The case of $x_{PQ}^\gl$. } The required conditions on $x_{PQ}^\gl\in\R^2$ are 
that
\[
H_1(x_{PQ}^\gl)=0 \ \ \AND  \ \ v_0\cdot x_{PQ}^\gl=\gl. 
\] 

If we set $x_{PQ}^\gl=(s,t)$, then the conditions above are equivalent to that
\[
s=h_1(t) \ \ \AND \ \ v_{0}\cdot(s,t)=\gl,
\]
which can be stated as 
\[
s=h_1(t) \ \ \AND \ \ v_{0}\cdot(h_1(t),t)=\gl. 
\]
If $t\in\R$ satisfies the above conditions, then $t>0$. Indeed, if $t\leq 0$ were the case, then we would have $v_{0}\cdot(h_1(t),t)\leq 0$, which is impossible. 
Hence, the unique existence of $x_{PQ}^\gl$ with the required properties is 
equivalent to the unique existence of the solution $t>0$ of 
$v_{0}\cdot(h_1(t),t)=\gl$, which can be written as  $f_1(t)=\gl$.   By Lemma \ref{4 func}, the solution of $f_1(t)=\gl$ is given uniquely by $t=f_1^{-1}(\gl)$. 
Hence, $x_{PQ}^\gl$ is given uniquely as the point $(h_1\circ f_1^{-1}(\gl),f_1^{-1}(\gl))$. That is, $x_{PQ}^\gl=G_1\circ f_1^{-1}(\gl)$. 

{\sc The case of $x_{PR}^\gl$. } Similarly to the above case, 
the problem of finding $x_{PR}^\gl$ is stated as that of finding 
a solution $t>0$ of $f_2(t)=\gl$, which can be solved uniquely 
by $t=f_2^{-1}(\gl)$ and the point $x_{PR}^\gl$ is given by $G_2\circ f_2^{-1}(\gl)$.

{\sc The case of $x_{QS}^\gl$. } Set $\gl_1:= v_1\cdot x_{PQ}^\gl$ and $t=f_1^{-1}(\gl)$.  Since $t>0$, $x_{PQ}^\gl=G_1(t)$, and the function $g_1\mid s\mapsto v_1\cdot G_1(s)$ is strictly increasing on $[0,\infty)$ by Lemma \ref{4 func}, we find that 
$\gl_1=g_1(t)>0$.  
The conditions on $x_{QS}^\gl$ are now stated as 
\[
H_2(x_{QS}^\gl)=0, \ \  v_1\cdot x_{QS}^\gl=\gl_1, \ \ \AND \ \ 
x_{QS}^\gl\in (-\infty,0)^2. 
\]
The first condition above, together with the last one, can be stated as $x_{QS}^\gl=(s,h_2(s))$ for some $s<0$.  Hence, 
finding a point $x_{QS}^\gl$ with the properties above is equivalent to finding 
a solution $s<0$ of $p_1(s)=\gl_1$, which can be solved uniquely as 
$s=p_1^{-1}(\gl_1)$. Thus, there is a unique point $x_{QS}^\gl$ having the required 
properties. Moreover, we find the formula for $x_{QS}^\gl$:
\beq \label{x_QS}
x_{QS}^\gl=G_2\circ p_1^{-1}(\gl_1)=G_2\circ p_1^{-1}\circ g_1\circ f_1^{-1}(\gl).
\eeq

{\sc The case of $x_{RS}^\gl$. } We skip the detail, but we point out that $x_{RS}^\gl$ is given by the formula
\beq\label{x_RS}
x_{RS}^\gl=G_1\circ p_2^{-1}\circ g_2\circ f_2^{-1}(\gl).
\eeq
The proof is now complete. 
\eproof 

Using formulas \erf{fgpq}, we can define $f_i,g_i,p_i,q_i$, $i=1,2$, on $\R$.  
The convexity property of functions $f_i, g_i, p_i, q_i$, $i=1,2$, are recorded 
in the following. 

\begin{lemma}\label{conv fgpq} The functions $f_i,\,-g_i,\,p_i,\,-q_i$, $i=1,2$, 
are convex on $\R$. 
\end{lemma}

\bproof The differentiation in the distributional sense twice yields 
\[\bald
&f_1''=v_{01}h_1''\geq 0, \quad
f_2''=v_{02}h_2''\geq 0,\quad
g_1''=v_{11}h_1''\leq 0,\quad
g_2''=v_{22}h_2''\leq 0,\quad
\\& p_1''=v_{12}h_2''\geq 0, \quad
p_2''=v_{21}h_1''\geq 0,\quad
q_1''=-\gam h_2''\leq 0,\quad
q_2''=-\gam h_1''\leq 0,
\eald \]
which assures the convexity of the functions \ 
$f_1, \, f_2, \, -g_1,\, -g_2,\, p_1,\, p_2,\, -q_1,\,-q_2$\ on\ $\R$.
\eproof

\section{Additional assumption}

It is not clear for the authors if conditions \erf{cond-H_i} 
are enough to build a smooth test function, by which one can 
prove the comparison principle for solutions of \erf{N-prob}.  
We here make another assumption on the functions $H_i$ in addition to 
\erf{cond-H_i}--\erf{zero-norm}, 
as described below.

With the notation above, we introduce an additional condition on the functions $H_i$
and the choice of $v_i,\, i=0,1,2$. 

\beq\tag{B3}\label{cond-slope} \left\{\ \begin{minipage}{0.8\textwidth}
There exists a constant $0<\gam<1$ such that for each $\gl>0$, the slope of the line passing through $x_{QS}^\gl$ and $x_{RS}^\gl$ ranges in the interval $[-1/\gam,-\gam]$. 
\end{minipage}\right.
\eeq

In our presentation under the assumption \erf{cond-slope}, 
$m_\gam$ and $n_\gam$ are those with $\gam$ given by \erf{cond-slope}. 
 
We give here three cases where \erf{cond-slope} holds.

\begin{example} \label{ex1}Assume that \erf{cond-H_i}--\erf{zero-norm} hold and that the functions $h_i$ 
are positively homogeneous of degree one, that is, for some constants $\ga_i^0
\in [0,\ga_i]$, 
\[
h_i(t)= \bcases
\ga_i^0 t \ & \IF \ t\leq 0,\\
\ga_i t & \IF \ t\geq 0.
\ecases
\] 
Then \erf{cond-slope} holds for some choice of $v_0,v_1,v_2\in\R^2$ satisfying \erf{cond-v_i},\erf{cond-v_1}, and \erf{cond-v_2}.
\end{example}

\bproof Consider the points 
\[
\tilde x_{PQ}=(\ga_1,1),\quad \tilde x_{PR}=(1,\ga_2), \quad
\tilde x_{QS}=(-1,-\ga_2^0),\quad \tilde x_{RS}=(-\ga_1^0,-1),
\]
which are on the half lines $x_1=\ga_1x_2\ (x_2>0)$, 
$x_2=\ga_2 x_1\ (x_1>0)$, $x_2=\ga_2^0 x_1 \ (x_1>0)$, and $x_1=\ga_1^0 x_2\ (x_2<0)$, respectively. The lines passing through the two points $\tilde x_{PQ}$ 
and $\tilde x_{PR}$, $\tilde x_{PQ}$ and $\tilde x_{QS}$, and $\tilde x_{PR}$ and $\tilde x_{RS}$, are given, respectively, by 
\[
(1-\ga_2)x_1+(1-\ga_1)x_2=1-\ga_1\ga_2, \quad -(1+\ga_2^0)x_1+(1+\ga_1)x_2=1-\ga_1\ga_2^0,
\]
and
\[
(1+\ga_2)x_1-(1+\ga_1^0)x_2=1-\ga_1^0\ga_2. 
\]
We set
\[
v_0=\Big(\fr{1-\ga_2}{1-\ga_1\ga_2},\fr{1-\ga_1}{1-\ga_1\ga_2}\Big), 
\quad
v_1=\Big(-\fr{1+\ga_2^0}{1-\ga_1\ga_2^0},\fr{1+\ga_1}{1-\ga_1\ga_2^0}\Big), 
\]
and
\[
v_2=\Big(\fr{1+\ga_2}{1-\ga_1^0\ga_2}, -\fr{1+\ga_1^0}{1-\ga_1^0\ga_2}\Big). 
\]
It is now easy to check that \erf{cond-v_i}, \erf{cond-v_1}, and \erf{cond-v_2} hold
and that for each $\gl>0$, 
\[
x_{PQ}^\gl=(\ga_1\gl,\gl),\quad x_{PR}^\gl=(\gl,\ga_2\gl),
\quad x_{QS}^\gl=(-\gl,-\ga_2^0\gl), \quad x_{RS}^\gl=(-\ga_1^0\gl,-\gl). 
\]
Hence, the slope $s$ of the line passing through $x_{QS}^\gl$ and $x_{RS}^\gl$ is given by 
\[
s=-\fr{1-\ga_2^0}{1-\ga_1^0},
\]
and, accordingly, if 
\[
0<\gam\leq \min\Big\{\fr{1-\ga_2^0}{1-\ga_1^0},\fr{1-\ga_1^0}{1-\ga_2^0}\Big\},
\]
then $s\in[-1/\gam,-\gam]$. Thus, \erf{cond-slope} holds with the current choice 
of $v_i$, $i=0,1,2$. 
\eproof

A small perturbation of the $h_i$ in the example above keeps the validity of \erf{cond-slope} 
as stated in the next example. 

\begin{example}Assume that \erf{cond-H_i}--\erf{zero-norm} hold. For $i=1,2$, set 
\[
\ga_i^-=\einf_{t\in\R}h_i'(t), \quad \ga_i^{0-}=\esup_{t\leq 0}h_i'(t) \ \ \AND \ \  
\ga_i^{0+}=\einf_{t\geq 0}h_i'(t). 
\]
Then there is a constant $\ep>0$ such that if, for $i=1,2$,
\[
\max\{\ga_i-\ga_i^{0+}, \ga_i^{0-}-\ga_i^-\}\leq \ep, 
\] then \erf{cond-slope} holds. 
\end{example}

\bproof  We start by noting that if $\ga_i^-=\ga_i^{0-}$ and $\ga_i^{0+}=\ga_i$ 
for $i=1,2$, then the situation is the same as Example \ref{ex1}, with 
$\ga_i^0=\ga_i^{0-}$, $i=1,2$. 

As in Example \ref{ex1}, we set \[
v_0=\Big(\fr{1-\ga_2}{1-\ga_1\ga_2},\fr{1-\ga_1}{1-\ga_1\ga_2}\Big), 
\quad
v_1=\Big(-\fr{1+\ga_2^{0-}}{1-\ga_1\ga_2^{0-}},\fr{1+\ga_1}{1-\ga_1\ga_2^{0-}}\Big), 
\]
and
\[
v_2=\Big(\fr{1+\ga_2}{1-\ga_1^{0-}\ga_2}, -\fr{1+\ga_1^{0-}}{1-\ga_1^{0-}\ga_2}\Big),
\]
and note that \erf{cond-v_i}, \erf{cond-v_1}, and \erf{cond-v_2} hold.

Fix any $\gl>0$. Writing $x_{PQ}^\gl=(x_1,x_2)$, we have 
\[
(1-\ga_2)x_1+(1-\ga_1)x_2=(1-\ga_1\ga_2)\gl \ \ \ \AND \ \ \ x_1=h_1(x_2).
\]
Since $\ga_1^{0+}x_2\leq h_1(x_2)\leq \ga_1 x_2$, there is $\gb_1^+\in[\ga_1^{0+},\ga_1]$ such that $h_1(x_2)=\gb_1^+x_2$, and then,
\[
x_1=\fr{\gb_1^+(1-\ga_1\ga_2)\gl}{\gb_1^+(1-\ga_2)+1-\ga_1} \ \ \AND \ \ x_2=\fr{(1-\ga_1\ga_2)\gl}{\gb_1^+(1-\ga_2)+1-\ga_1}. 
\]
Similarly, if $x_{PR}^\gl=(x_1,x_2)$, then we have for some $\gb_2^+\in[\ga_2^{0+},\ga_2]$,
\[h_2(x_1)=\gb_2^+ x_1, \ \ \ 
x_1=\fr{(1-\ga_1\ga_2)\gl}{1-\ga_2+\gb_2^+(1-\ga_1)} \ \  \AND \ \ 
x_2=\fr{\gb_2^+(1-\ga_1\ga_2)\gl}{1-\ga_2+\gb_2^+ (1-\ga_1)}. 
\]
Setting 
\[
\gth_1=\fr{1-\ga_1\ga_2}{1-\ga_1+\gb_1^+(1-\ga_2)} \ \ \AND \ \ \gth_2
=\fr{1-\ga_1\ga_2}{1-\ga_2+\gb_2^+(1-\ga_1)},
\]
we have 
\[
x_{PQ}^\gl=(\gb_1^+\gth_1\gl,\gth_1\gl) \ \  \AND \ \ x_{PR}^\gl=(\gth_2 \gl,\gb_2^+\gth_2 \gl). 
\]

If we set 
\[
\mu:=\left(-(1+\ga_2^{0-}), 1+\ga_1\right)\cdot x_{PQ}^\gl \ \ \AND \ \ x_{QS}^\gl=:(x_1,x_2),
\]
then 
\[
-(1+\ga_2^{0-})x_1+(1+\ga_1)x_2=\mu \ \ \AND \ \ x_2=h_2(x_1). 
\]
Noting that $h_2(x_1)=\gb_2^- x_1$ for some $\gb_2^-\in[\ga_2^-,\ga_2^{0-}]$, we 
find that 
\[\begin{gathered}
\mu=\left(1+\ga_1 -\gb_1^+(1+\ga_2^{0-})\right)\gth_1\gl,
\\x_1=-\fr{\mu}{1+\ga_2^{0-}-\gb_2^-(1+\ga_1)} \ \ \AND 
\ \ x_2=-\fr{\gb_2^-\mu}{1+\ga_2^{0-}-\gb_2^-(1+\ga_1)}. 
\end{gathered}
\]
Similarly, we find that if we set 
\[
\nu=\left(1+\ga_2-\gb_2^+(1+\ga_1^{0-}\right)\gth_2\gl,
\]
then for some $\gb_1^-\in[\ga_1^-,\ga_1^{0-}]$,
\[
x_{RS}^\gl=\Big(-\fr{\gb_1^-\nu}{1+\ga_1^{0-}-\gb_1^-(1+\ga_2)},-\fr{\nu}{
1+\ga_1^{0-}-\gb_1^-(1+\ga_2)}\Big). 
\]

We note that, as $\ga_i^- \to \ga_i^{0-}$ and 
$\ga_i^{0+}\to\ga_i$, 
\[\bald
&\Big(-\fr{\mu}{1+\ga_2^{0-}-\gb_2^-(1+\ga_1)},-\fr{\gb_2^-\mu}{1+\ga_2^{0-}-\gb_2^-(1+\ga_1)}\Big)\ \to \ (-\gl,-\ga_2^{0-}\gl),
\\ &\Big(-\fr{\gb_1^-\nu}{1+\ga_1^{0-}-\gb_1^-(1+\ga_2)},-\fr{\nu}{
1+\ga_1^{0-}-\gb_1^-(1+\ga_2)}\Big) \ \to \ (-\ga_1^{0-}\gl, -\gl). 
\eald
\] 
We choose constants $\kappa_i$, $i=1,2$, so that 
\[
0<\kappa_1<\min_{i=1,2}\ga_i^{0-}<\kappa_2<1.
\]
By the convergence above, we may select $\ep>0$ small enough so that if
\beq\label{ex2.1}
\max_{i=1,2}\max\{\ga_i^{0-}-\ga_i^-,\ga_i-\ga_i^{0+}\}\leq \ep,
\eeq
then 
\[\bald
&-\fr{\mu}{1+\ga_2^{0-}-\gb_2^-(1+\ga_1)}, -\fr{\nu}{
1+\ga_1^{0-}-\gb_1^-(1+\ga_2)}\in[-2\gl,\, -\kappa_2\gl],
\\& 
-\fr{\gb_2^-\mu}{1+\ga_2^{0-}-\gb_2^-(1+\ga_1)}, 
-\fr{\gb_1^-\nu}{1+\ga_1^{0-}-\gb_1^-(1+\ga_2)}\in [-\kappa_1\gl,\, 0].
\eald
\]
Thus, if \erf{ex2.1} is valid, then the slope $s$ of the line passing through 
$x_{QS}^\gl$ and $x_{RS}^\gl$ satisfies
\[
s\in\Big[-\fr{2}{\kappa_2-\kappa_1},\,-\fr{\kappa_2-\kappa_1}{2}\Big],
\]
which shows that \erf{cond-slope} holds if \erf{ex2.1} is satisfied. 
\eproof

\begin{example} Assume that \erf{cond-H_i}--\erf{zero-norm} hold and moreover that either 
$h_1(t)=0$ for all $t\leq 0$ or $h_2(t)=0$ for all $t\leq 0$. Then \erf{cond-slope} holds.
\end{example}

\bproof  By symmetry, we may only consider the case when $h_2(t)=0$ for all 
$t\leq 0$. Fix $v_0=(1,1)$ and $v_2=(1,-1)$, so that \erf{cond-v_i} for $v_0$ and $v_2$
and \erf{cond-v_2} hold.  Set 
\[
\ga_1^{0}=\esup_{t\leq 0}h_1'(t),
\]
fix a constant $\ep_0>0$, and set  
\[
\gb_0=\ep_0+\fr{\ga_1^0}{1-\ga_1^0}.
\]
Set 
\[
\gb=\ga_1+\gb_0 (1+\ga_1) \ \ \AND \ \ v_1=(-1,\gb),
\]
and note that $\gb\geq \gb_0\geq \ep_0>0$ 
and that \erf{cond-v_i} for $v_1$ and \erf{cond-v_1} hold. 

Fix any $\gl>0$ and let $x_{PQ}^\gl,x_{QS}^\gl, x_{PR}^\gl, x_{RS}^\gl$ denote the points 
described above (see also Proposition \ref{four-points}). By solving 
\[
v_0\cdot x=\gl,\quad x_2=h_2(x_1),
\]
and writing $\xi_1^\gl$ for the unique solution $t$ of $t+h_2(t)=\gl$, 
we have
\[
x_{PR}^\gl=(\xi_1^\gl,h_2(\xi_1^\gl)) \ \ \AND \ \ \fr \gl {1+\ga_2}<\xi_1^\gl\leq \gl.
\]
Similarly, we have 
\[
x_{PQ}^\gl=(h_1(\xi_2^\gl),\xi_2^\gl) \ \ \AND \ \ \fr{\gl}{1+\ga_1}<\xi_2^\gl\leq \gl,
\]
where $\xi_2^\gl$ is the unique solution $t$ of $t+h_1(t)=\gl$. Also, we have 
\[
x_{QS}^\gl=(\eta_1^\gl,0) \ \ \AND \ \ x_{RS}^\gl=(h_1(\eta_2^\gl),\eta_2^\gl),
\]
where $\eta_1^\gl:=-v_1\cdot x_{PQ}^\gl$ 
and $\eta_2^\gl$ denotes the unique solution 
$t$ of $t-h_1(t)=-v_2\cdot x_{PR}^\gl$. 

Noting that $\eta_i^\gl<0$, $i=1,2$, 
\[\bald
\eta_1^\gl&=-v_1\cdot x_{PR}^\gl =h_1(\xi_2^\gl)-\gb\xi_2^\gl,
\\ \eta_2^\gl&=h_1(\eta_2^\gl)-v_2\cdot x_{PR}^\gl=h_1(\eta_2^\gl)-\xi_1^\gl+h_2(\xi_1^\gl),
\eald
\]
we compute that 
\[
\eta_1^\gl=h_1(\xi_2^\gl)-\gb \xi_2^\gl 
\left\{
\bald
&\leq (\ga_1-\gb)\xi_2^\gl=-\gb_0 (1+\ga_1) \xi_2^\gl\leq -\gb_0\gl,
\\& \geq -\gb \xi_2^\gl\geq -\gb\gl,
\eald
\right.
\]
and
\[
\eta_2^\gl\leq -\xi_1^\gl+h_2(\xi_1^\gl) \leq -(1-\ga_2)\xi_1^\gl< -\fr{(1-\ga_2)\gl}{1+\ga_2}. 
\]
Also, we have 
\[
\eta_2^\gl\geq \ga_1^0\eta_2^\gl-\xi_1^\gl\geq \ga_1^0\eta_2^\gl -\gl,
\]
and hence, 
\[
\eta_2^\gl\geq -\fr{\gl}{1-\ga_1^0}. 
\]
Furthermore, we note that
\[
\eta_1^\gl-h_1(\eta_2^\gl) >\eta_1^\gl\geq -\gb \gl,
\]
and
\[
\eta_1^\gl-h_1(\eta_2^\gl)\leq \eta_1^\gl-\ga_1^0\eta_2^\gl 
\leq-\gb_0\gl+\ga_1^0\fr{\gl}{1-\ga_1^0}
=(\ga_1^0-\gb_0(1-\ga_1^0))\fr{\gl}{1-\ga_1^0}=-\ep_0 \gl.
\]

If we denote by $s$ the slope of the line passing through $x_{QS}^\gl$ and $x_{RS}^\gl$, then
\[
s=\fr{-\eta_2^\gl}{\eta_1^\gl-h_1(\eta_2^\gl)}
\left\{\bald
&\leq \fr{-\eta_2^\gl}{-\gl\gb}\leq -\fr{1-\ga_2}{(1+\ga_2)\gb}, 
\\&\geq  \fr{-\eta_2^\gl}{-\ep_0\gl}\geq -\fr{1}{(1-\ga_1^0)\ep_0}.
\eald
\right.
\]
It is immediate to see that \erf{cond-slope} is satisfied with any positive 
constant $\gam$ such that 
\[
\gam\leq \min\Big\{\fr{1-\ga_2}{(1+\ga_2)\gb},(1-\ga_1^0)\ep_0\Big\}. \qedhere
\]
\eproof

 Fix any point $y=G_2(t)$, with $t<0$, on the curve $H_2=0$. 
Let $s_0\geq 0$ be the slope of the line passing through the origin and the point $y$. 
Let $G_1(s)$ be an arbitrary point on the curve $H_1=0$ and let $\gs(s)$ denote the slope of the line passing through the points $y$ and $G_1(s)$, which makes sense as far as $t\not=h_1(s)$. Since $h_1$ is nondecreasing, 
$h_1(-\infty):=\lim_{s\to-\infty}h_1(s)\in[-\infty,0]$ makes sense. Since $h_1$ is convex and $h_1(0)=0$, if $h_1(-\infty)<t$, then there is a unique $\tau\in(-\infty,0)$ such that 
$h_1(\tau)=t$. Similarly, if $h_1(-\infty)=t$, then there is a unique 
$\tau\in[-\infty,0)$ such that $h_1(\tau)=t$ and $h_1(s)>t$ for every $s>\tau$. 
If $h_1(-\infty)>t$, then we set $\tau=-\infty$. 
We record the following simple observation. 
  
\begin{lemma} \label{func-gs}With the notation introduced above, the function $\gs$ 
is increasing and continuous on $(\tau,0]$, $\gs(0)=s_0$, and $\lim_{s\to \tau+} \gs(s)=-\infty$. Furthermore, if $h_1(-\infty)<t$, then $\gs(s)>0$ on 
$(-\infty,\tau)$ and if $h_1(-\infty)=t$ and $\tau>-\infty$, then for any $s\in(-\infty,\tau]$,  $h_1(s)=t$ and $\gs(s)$ does not make sense. 
\end{lemma}

\bproof By the definition of $\gs$, 
if $h_1(s)\not=t$, then 
\[
\gs(s)=\fr{h_2(t)-s}{t-h_1(s)}. 
\]
It is obvious that $\gs$ is continuous on $(\tau,0]$. 
By the choice of $s_0$, we have $\gs(0)=s_0$. 

Note that since $0\leq h_1'<1$ and $0\leq h_2'<1$ a.e. and $h_1(0)=h_2(0)=0$, we have  
$\tau<t$ and $t<h_2(t)$.   Note also by the definition of $\tau$ that $t<h_1(s)$ for any $s\in(\tau,0]$.

To show the monotonicity of $\gs$ on $(\tau,0]$, we argue slightly in an informal way, 
that is, we assume that $h_1\in C^2(\tau,0]$.  
We compute that for any $s\in(\tau,0]$,
\[
\gs'(s)=\fr{-(t-h_1(s))+h_1'(s)(h_2(t)-s)}{(t-h_1(s))^2},
\]
and set
\[
\nu(s)=-(t-h_1(s))+h_1'(s)(h_2(t)-s).
\]
It follows that 
\[
\nu'(s)=h_1''(s)(h_2(t)-s).
\]
Note that for any $s\in(\tau,h_2(t)]$,
\[
\nu(s)\geq h_1(s)-t>0.
\]
Observe that $\nu'(s)\leq 0$ for all $s\in[h_2(t),0]$, which implies that 
$\nu$ is nonincreasing on $[h_2(t),0]$, and that 
\[
\nu(0)=-t+h_1'(0)h_2(t)=h_1'(0)(h_2(t)-t)+t(h_1'(0)-1)>h_1'(0)(h_2(t)-t)\geq 0.
\]
Hence, we have $\nu(s)>0$ and $\gs'(s)>0$ on $(\tau,0]$.  
Thus, $\gs$ is strictly increasing on $(\tau,0]$. The argument above can be justified 
with a bit more elaboration, which we leave to the interested reader. 

If $\tau>-\infty$, then it is obvious that
\[
\lim_{s\to\tau+}|\gs(s)|=\infty,
\]
and moreover, by the monotonicity of $\gs$, 
\[
\lim_{s\to\tau+}\gs(s)=-\infty. 
\]

If $h_1(-\infty)<t$, then we have $t=h_1(\tau)$ and
and for any $s\in(-\infty,\tau)$,
\[
t-h_1(s)>0 \ \ \AND \ \  h_2(t)-s >h_2(t)-\tau,
\]
which yields \ $\gs(s)>0$.

If $h_1(-\infty)=t$ and $\tau>-\infty$, then $h_1(s)=t$ for all $s\in(-\infty,\tau]$
due to the convexity and monotonicity of $h_1$, and hence, $\gs(s)$ does not make sense for $s\in(-\infty,\tau]$.  
\eproof

\begin{lemma} \label{comp-slope} Let $\gl>0$ and let $t,s\in\R$ be such that $x_{QS}^\gl=G_2(t)$ 
and $x_{RS}^\gl=G_1(s)$. Set $\mu=m_\gam\cdot x_{QS}^\gl$ and $\nu=n_\gam\cdot x_{RS}^\gl$. Let $y_1=G_2(r_1)$ and $y_2=G_1(r_2)$ be the  unique points, respectively,  where the line $n_\gam\cdot x=\nu$ intersects 
the curve $H_1=0$ and where the line $m_\gam\cdot x=\mu$ intersects 
the curve $H_2=0$. Then, we have \ $r_1\leq t$ and $r_2\leq s$. 
\end{lemma}

The unique existence of $y_i$, $i=1,2$, can be proved similarly to that of 
$x_{PQ}^\gl,\,x_{PR}^\gl$, etc.

\bproof  We first remark that the lines $m_\gam\cdot x=\mu$ 
and $n_\gam\cdot x=\nu$ have slopes $-1/\gam$ and $-\gam$, respectively. 
Since $x_{QS}^\gl=G_2(t), x_{RS}^\gl=G_1(s)\in(-\infty,0)^2$, we have 
$t,s<0$ and $\mu, \nu>0$. Since $m_\gam\cdot G_2(r_1)=\mu$ 
and $n_\gam\cdot G_1(r_2)=\nu$, it is obvious that $r_1,r_2<0$. 
As in Lemma \ref{func-gs}, let  $\gs(\xi)$ denote the slope of the line passing through the points $x_{QS}^\gl$ and $G_1(\xi)$ as far as the ``slope" makes sense. 
We have $\gs(r_2)=-1/\gam$ and, by \erf{cond-slope}, 
$\gs(s)\in[-1/\gam,-\gam]$. The monotonicity of $\gs$ by Lemma 
\ref{func-gs} assures that $r_2\leq s$. A parallel argument to the above yields 
$r_1\leq t$.  
\eproof

For later convenience, we remark that, with the notation of Lemma \ref{comp-slope} and \erf{fgpq}, 
we have 
\beq\label{gl to mu}\left\{\ \bald
&t=f_1^{-1}(\gl), \quad x_{QS}^\gl=G_2\circ p_1^{-1}\circ g_1\circ f_1^{-1}(\gl), \quad
x_{RS}^\gl=G_1\circ p_2^{-1}\circ g_2\circ f_2^{-1}(\gl),
\\&  \mu=q_1\circ  p_1^{-1}\circ g_1\circ f_1^{-1}(\gl),\quad \nu
=q_2\circ p_2^{-1}\circ g_2\circ f_2^{-1}(\gl). 
\eald \right. 
\eeq

\section{Convex function $\psi$}

Our first step of constructing the test function $\phi_\ep$  in Theorem \ref{exist-phi}
is to build a nonnegative convex function $\psi$, the square of which is a good approximation of the convex conjugate function of $\phi_\ep$.  

The idea of building $\psi$ is to let $\psi$ be a convex function that has a  linear 
growth and that for each $\gl>0$, the level set $\{\psi=\gl\}$ includes the line segments $[x_{PQ}^\gl, x_{PR}^\gl]$, 
$[x_{QS}^\gl,x_{PQ}^\gl]$, and $[x_{RS}^\gl,x_{PR}^\gl]$ (see Proposition \ref{four-points} for the definition of $x_{PQ}^\gl$, $x_{PR}^\gl$, etc.). 
Here, $[x,y]$ denotes the set $\{(1-t)x+ty\mid 0\leq t\leq 1\}$ for $x,y\in\R^2$. 
For notational 
convenience, we write 
\[
I_P^\gl=[x_{PQ}^\gl, x_{PR}^\gl],\quad I_Q^\gl= [x_{QS}^\gl,x_{PQ}^\gl], \ \ \AND \ \ 
I_R^\gl=[x_{RS}^\gl,x_{PR}^\gl].
\]

\begin{lemma}\label{I^gl} We have 
\[
P\stm\{0\}=\bigcup_{\gl>0}I_P^\gl,\quad 
Q\stm\{0\}=\bigcup_{\gl>0}I_Q^\gl, \ \ \AND \ \ R\stm\{0\}
=\bigcup_{\gl>0}I_R^\gl.
\]
\end{lemma}

\bproof Note first that 
\[
\pl P=(\{x\in \R^2\mid H_1(x)=0\} \cup\{x\in\R^2 \mid H_2(x)=0\}) \cap\{x\in \R^2\mid v_0\cdot x\geq 0\}. 
\] 
Fix any $\gl>0$, note that 
that 
\[
\{(1-t)x_{PQ}^\gl+t x_{PR}^\gl\mid t\in\R\}=\{x\in\R^2\mid v_0\cdot x=\gl\}, 
\]
and, in particular,
\[
I_P^\gl\subset \{x\in\R^2\mid v_0\cdot x=\gl\}. 
\]
Recall by Proposition \ref{four-points} that the points $x_{PQ}^\gl$  and $x_{PR}^\gl$ are the unique points where the line 
$v_0\cdot x=\gl$ meets, respectively, the curves $H_1=0$ and $H_2=0$. 
Set 
\[
I=\{t\in\R\mid (1-t)x_{PQ}^\gl+t x_{PR}^\gl\in P\},
\]
and note that $I$ is a compact subset of $\R$ (since the intersection of the line $v_0\cdot x=\gl$ and $P$ is compact). In particular,  we have 
$0,1\in I\subset [\min I,\max I]$
The uniqueness of the points $x_{PQ}^\gl$  and $x_{PR}^\gl$ shows that if $a\in \pl I$, then $(1-a)x_{PQ}^\gl+ax_{PR}^\gl\in \pl P$ and hence either $a=0$ or $a=1$.   
It is clear that the line $v_0\cdot x=\gl$ intersects the interior of $P$, which implies 
that $I$ has a nonempty interior. 
We now deduce that $I=[0,1]$, which reads
\[
I_P^\gl=P\cap\{x\in \R^2\mid v_0\cdot x=\gl\}. 
\]
Thus, we conclude that 
\[
P\stm \{0\}=P\cap\{x\in\R^2\mid v_0\cdot x>0\} 
=\bigcup_{\gl>0}P\cap \{x\in \R^2\mid v_0\cdot x=\gl\}
=\bigcup_{\gl>0}I_P^\gl.
\]

Let $\gl>0$ and note by Proposition \ref{four-points} (and its proof) 
that $x_{PQ}^\gl$ and $x_{QS}^\gl\in (-\infty,0)^2$ are 
the unique points where the line $v_1\cdot x=\gl_1$, with 
$\gl_1=g_1\circ f_1^{-1}(\gl)$, meets the curves $H_1=0$ and $H_2=0$, 
respectively.   Note also that $x_{PR}^\gl$ and $x_{RS}^\gl$ are the unique points where the line $v_2\cdot x=\gl_2$, with $\gl_2=g_2\circ f_2^{-1}(\gl)$, 
meets the curves $H_2=0$ and $H_1=0$, respectively. 
Arguing as above, we find that 
\[
I_Q^\gl=Q\cap\{x\in \R^2\mid v_1\cdot x=\gl_1\} \ \ \AND \ \ 
I_R^\gl=R\cap\{x\mid v_2\cdot x=\gl_2\}. 
\]
Since both $g_i\circ f_i^{-1}$, $i=1,2$, map $(0,\infty)$ onto itself, we conclude that 
\[
Q\stm\{0\}=\bigcup_{\gl>0}I_Q^\gl \ \ \AND \ \ R\stm\{0\} 
=\bigcup_{\gl>0}I_R^\gl. \qedhere
\]
\eproof 

Define the functions $\psi_0$, $\psi_1$, and $\psi_2$ on $\R^2$ by 
\beq\label{d-psi_i}
\psi_0(x)=(v_0\cdot x)_+, \quad 
\psi_1(x)=f_1\circ g_1^{-1}((v_1\cdot x)_+), \quad \AND \quad
\psi_2(x)=f_2\circ g_2^{-1}((v_2\cdot x)_+),
\eeq
where $a_+$ denotes the positive part, $\max\{a,0\}$, of $a\in\R$. 
It is clear that $\psi_0$ is convex on $\R^2$. By Lemma \ref{4 func}, we find that 
the functions $f_i\circ g_i^{-1}$, $i=1,2$, are uniformly increasing and Lipschitz continuous on $[0,\infty)$. In particular, for some constants $0<\gd<C$, 
\beq\label{psi_i l-growth}
\gd(v_i\cdot x)_+\leq \psi_i(x)\leq C(v_i\cdot x)_+ \ \ \FORALL x\in\R^2.
\eeq
By Lemma \ref{conv fgpq}, $f_i$, $i=1,2$, are convex and $g_i$, $i=1,2$, are 
concave. Since $f_i, g_i$, $i=1,2$, are increasing, it follows that 
$g_i^{-1}$, $i=1,2$, are convex on $[0,\infty)$ and moreover, 
$f_i\circ g_i^{-1}$, $i=1,2$, are convex on $[0,\infty)$. Noting that for $i=1,2$, 
$x\mapsto (v_i\cdot x)_+$ is convex on $\R^2$, we find that 
$\psi_i$, $i=1,2$, are convex on $\R^2$.

Let $\gl>0$. We next observe that 
\beq\label{psi_012}
\psi_0(x)=\gl \ \ \ON I_P^\gl, 
\quad 
\psi_1(x)=\gl \ \ \ON I_Q^\gl,
 \ \ \AND \ \ 
\psi_2(x)=\gl \ \ \ON I_R^\gl.
\eeq
By definition, we have $v_0\cdot x_{PQ}^\gl=v_0\cdot x_{PR}^\gl=\gl$. 
In particular, if $x\in I_P$, we have $v_0\cdot x=\gl$ and 
$\psi_0(x)=\gl$. 
If $x_{PQ}^\gl=G_1(t)$ and $x_{PR}^\gl=G_2(s)$ for some $t,s\in\R$, then 
$t,s>0$ and $\gl=f_1(t)=f_2(s)$. Setting $\gl_1:=v_1\cdot x_{PQ}^\gl
=v_1\cdot G_1(t)$ and 
$\gl_2:=v_2\cdot x_{PR}^\gl=v_2\cdot G_2(s)$, by the definition of $g_i$, we have $g_1(t)=\gl_1$ and $g_2(s)=\gl_2$, which can be stated as $\gl_1=g_1^{-1}(t)$ 
and $\gl_2=g_2^{-1}(s)$. Hence, if $x\in I_Q^\gl$, then 
$v_1\cdot x=\gl_1$ and $\psi_1(x)=f_1\circ g_1^{-1}(\gl_1)=f_1(t)=\gl$. 
Similarly, if $x\in I_{R}^\gl$, then $\psi_2(x)
=f_2\circ g_2^{-1}(\gl_2)=f_2(s)=\gl$.

Set 
\[
\tilde \psi(x)=\max_{i=1,2,3}\psi_i(x) 
\ \ \FOR x\in\R^2.
\]
Obviously, $\tilde \psi$ is a nonnegative convex function on $\R^2$. 
We define $\psi\mid \R^2\to \R$ by
\beq\label{def-psi}
\psi(x)=\max\{\phi(x)\mid \phi \text{ convex on }\R^2,\ \phi\leq \tilde \psi \ON P\cup Q\cup R\}.
\eeq
It is clear that $\psi$ is a convex function on $\R^2$ and that 
$\tilde\psi(x)\leq \psi(x)$ for all $x\in\R^2$ and $\psi(x)= \tilde \psi(x)$ for all $x\in P\cup Q\cup R$.

\begin{lemma} \label{psi=psi_012} We have 
\[
\psi=\psi_0 \ \ON \ P,\quad \psi=\psi_1 \ \ON \ Q,\ \ \AND \ \  \psi=\psi_2 \ \ON \ R.
\]
\end{lemma}

\bproof Note that $\psi=\tilde \psi$ on $P\cup Q\cup R$. 
Let $\gl>0$ and let $y_1^\gl, \, y_2^\gl$ be the unique points where the line $v_0\cdot x=\gl$ intersects the lines $x_1=0$ (i.e., the $x_2$-axis) 
and $x_2=0$ (the $x_1$-axis), respectively. 
Consider the convex functions $\psi_i$, $i=0,1,2$, 
on the line segment $I_0:=[y_1^\gl,y_2^\gl]$. By the definition of $\psi_0$, we have 
$\psi_0(x)=\gl$ on $I_0$. 
Since $v_{22}<0$, we have $v_2\cdot y_1^\gl<0$ and $\psi_2(y_1^\gl)=
f_2\circ g_2^{-1}(0)=0$. Similarly, we have $\psi_1(y_2^\gl)=0$. 
Note that $I_P^\gl\subset I_0$ and that $\psi_1(x_{PQ}^\gl)=\gl$ 
and $\psi_2(x_{PR}^\gl)=\gl$. By the convexity of $\psi_i$, $i=0,1,2$, we infer that 
$\psi_1(x)\leq \psi_0(x)$ on $[x_{PQ}^\gl,y_2^\gl]$, 
$\psi_2(x)\leq \psi_0(x)$ on $[y_1^\gl,x_{PR}^\gl]$, and 
$\psi(x)=\tilde\psi(x)=\psi_0(x)$ on $I_P^\gl$.  

Set $\gl_1:=v_1\cdot x_{PQ}^\gl$ and note that $\gl_1=g_1\circ f_1^{-1}(\gl)>0$. 
As we have observed above, we have  $\psi_1(x_{PQ}^\gl)=\psi_0(x_{PQ}^\gl)=\gl\geq \psi_2(x_{PQ}^\gl)$.  Since $x_{QS}^\gl\in(-\infty,0)^2$, we have 
$v_0\cdot x_{QS}^\gl<0$ and $\psi_0(x_{QS}^\gl)=0$. 
We claim that $v_2\cdot x_{QS}^\gl<0$. Indeed, since $v_{22}<0$, if we define 
$g_2$ on $\R$ by $g_2(t)=v_2\cdot G_2(t)$, then $g_2$ is concave on $\R$ by Lemma \ref{conv fgpq}, 
$g_2(0)=0$, and $g_2(t)>0$ for all $t>0$ by Lemma \ref{4 func}. 
Hence, choosing $r<0$ so that $x_{QS}^\gl=G_2(r)$, we have 
\[
0=g_2(0)\geq \tfrac 12(g_2(r)+g_2(-r)) >\tfrac 12 g_2(r),
\]
which assures that $v_2\cdot x_{QS}^\gl=g_2(r)<0$. 
This shows that $\psi_2(x_{QS}^\gl)=0$. 
It now follows that $\psi_0(x)\leq \psi_1(x)$, 
$\psi_2(x)\leq \psi_1(x)$, 
and $\psi(x)=\psi_1(x)=\gl$ on  $I_Q^\gl$. 

An argument parallel to the above shows that 
$\max\{\psi_0(x),\psi_1(x)\}\leq \psi_2(x)=\psi(x)=\gl$ on 
$I_R^\gl$. 

Thus, noting that $\psi(0)=\psi_i(0)=0$ for all $i=0,1,2$, we conclude by Lemma \ref{I^gl} that 
\[
\psi=\psi_0 \ \ \ON P,\quad 
\psi=\psi_1 \ \ \ON Q,\ \ \AND \ \ 
\psi=\psi_2 \ \ \ON  R.  \qedhere
\]
\eproof 

\begin{lemma} \label{psi l-growth}There exist constants $0<\gd\leq C<\infty$ such that 
\[
\gd|x|\leq \psi(x)\leq C|x| \ \ \FORALL x\in\R^2.
\]
\end{lemma}

\bproof By \erf{psi_i l-growth}, there are constants $0<\gd_0<C_0<\infty$ such that 
\beq\label{psi_i l-growth2}
\gd_0(v_i\cdot x)_+\leq \psi_i(x)\leq C_0(v_i\cdot x)_+ \ \ \FORALL x\in \R^2,\, i=0,1,2. 
\eeq

We claim that for some constant $\gd_1>0$, 
\beq \label{lb-PQR}
\max_{i=0,1,2}v_i\cdot x\geq \gd_1|x| \ \ \FORALL x\in P\cup Q\cup R.
\eeq
Observe first that since $v_{01}>0$ and $v_{02}>0$, we have for any $x\in[0,\infty)^2$,
\[
v_0\cdot x\geq \min_{j=1,2}v_{0j}(x_1+x_2)\geq \min_{j=1,2}|v_{0j}||x|. 
\]  
Next,  since $v_{11}<0$ and $v_{12}>0$, 
we have for any $x\in (-\infty,0]\tim[0,\infty)$, 
\[
v_1\cdot x=|v_{11}||x_1|+|v_{12}||x_2|\geq \min_{j=1,2}|v_{1j}||x|.
\]
Similarly, we have for any $x\in[0,\infty)\tim(-\infty,0]$, 
\[
v_2\cdot x\geq \min_{j=1,2}|v_{2j}||x|. 
\]
Hence, we have 
\[
\max_{i=0,1,2}v_i\cdot x\geq \min_{i=0,1,2,\, j=1,2}|v_{ij}||x| 
\ \ \FORALL x\in [0,\infty)\tim \R \, \cup\,  \R\tim [0,\infty).
\] 

To conclude the validity of \erf{lb-PQR}, it remains to check the case where 
either $x\in Q\cap(-\infty,0]^2$ or $x\in R\cap(-\infty,0]^2$.     
Let $x\in Q\cap(-\infty,0]^2$. Since $v_{12}>0$ and $x_2\geq h_2(x_1)$, we have 
\[
v_1\cdot x \geq v_1\cdot G_2(x_1)=p_1(x_1).
\]
Hence, by the uniform monotonicity of $p_1$, a consequence of  Lemma \ref{4 func}, 
we find that for some constant $\ep_1>0$, 
\[
v_1\cdot x\geq p_1(x_1)\geq \ep_1|x_1|. 
\] 
By \erf{cond-H_i}, the function $h_2$ has the Lipschitz constant $\ga_2$ and $h_2(0)=0$. Accordingly, we have
\[
|x|^2\leq x_1^2+h_2(x_1)^2\leq (1+\ga_2^2)|x_1|^2,
\]
and hence, 
\[
v_1\cdot x\geq \fr{\ep_1}{|(1,\ga_2)|}|x|
\]
Similarly, if $x\in R\cap(-\infty,0]^2$, then we have
\[
v_2\cdot x \geq v_2\cdot G_1(x_2)=p_2(x_2)\geq \ep_2|x_2|,
\]
where $\ep_2$ is a positive constant, 
and moreover, 
\[
v_2\cdot x\geq \fr{\ep_2}{|(\ga_1,1)|}|x|.
\] 
Now, setting 
\[
\gd_1=\min\left\{\min_{i=0,1,2,\,j=1,2}|v_{ij}|, \fr{\ep_1}{|(1,\ga_2)|},\fr{\ep_2}{|(\ga_1,1)|}\right\}, 
\]
we find that \erf{lb-PQR} holds. Combining this with \erf{psi_i l-growth2}, we have 
\[
\psi(x)=\tilde\psi(x)\geq \gd_0\gd_1 |x| \ \ \FORALL x\in P\cup Q\cup R.
\]
Since the function $x\mapsto \gd_0\gd_1|x|$ is convex on $\R^2$, we conclude that 
\[
\psi(x)\geq \gd_0\gd_1|x| \ \ \FORALL x\in\R^2. 
\]

From \erf{psi_i l-growth2}, we find immediately that 
\[
\tilde\psi(x)\leq C_0\max_{i=0,1,2}|v_i||x| \ \ \FORALL x\in\R^2,
\]
which implies that
\[
\psi(x)\leq C_0\max_{i=0,1,2}|v_i||x| \ \ \FORALL x\in P\cup Q\cup R.
\]
Since $\psi$ is convex on $\R^2$ and $\pl (-\infty,0)^2\subset Q\cup R$, if $x\in(-\infty,0)^2$, then $x=(1-t)(|x_1|+|x_2|)(-e_1)+t(|x_1|+|x_2|)(-e_2)$, where $t=|x_2|/(|x_1|+|x_2|)$,  and 
\[\bald
\psi(x)&\leq (1-t)\psi((|x_1|+|x_2|)(-e_1))+t\psi((|x_1|+|x_2|)(-e_2))
\\&\leq C_0\max_{i=0,1,2}|v_i| (|x_1|+|x_2|)
\leq \sqrt 2 C_0\max_{i=0,1,2}|v_i| |x|. 
\eald\]
Thus, we find that 
\[
\psi(x)\leq \sqrt 2 C_0\max_{i=0,1,2}|v_i||x| \ \ \FORALL x\in\R^2. 
\]
This completes the proof. 
\eproof

We are interested in the sublevel sets $\{\psi\leq \gl\}$ of $\psi$, with $\gl>0$, and 
we introduce functions $\chi_i$, $i=1,2$, by
\[\bald
\chi_1(x)&=f_1\circ g_1^{-1}\circ p_1\circ q_1^{-1}((m_\gam\cdot x)_+), 
\\ \chi_2(x)&=f_2\circ g_2^{-1}\circ p_2\circ q_2^{-1}((n_\gam\cdot x)_+). 
\eald
\] 
For $i=1,2$, by Lemmas \ref{4 func} and \ref{conv fgpq}, the functions 
$f_i,\,g_i$ are increasing, $p_i,\,q_i$ are decreasing, $f_i,\,p_i$ are 
convex, and $g_i,\,q_i$ are concave on their respective domains of definition. Hence, 
$q_i^{-1}$ is concave on $[0,\infty)$, which implies that $p_i\circ q_i^{-1}$ is 
convex (and increasing) on $[0,\infty)$. Moreover, since $f_i,\,g_i^{-1}$ are increasing and convex on 
$[0,\infty)$, we find that $f_i\circ g_i^{-1}\circ p_i\circ q_i^{-1}$ is 
convex and increasing on $[0,\infty)$. Hence, the functions $\chi_i$, $i=1,2$, are convex on $\R^2$.

\begin{lemma} \label{chi-psi}We have $\,\chi_1=\psi\,$ on the curve $H_2=0$ in $S$, $\,\chi_2=\psi\,$ 
on the curve $H_1=0$ in $S$, and $\,\chi_i\leq \psi\,$ on $\R^2$ for $i=1,2$. 
\end{lemma}

Remark that the curve $H_2=0$ (resp., $H_1=0$) in $S$ is the set 
\[\begin{gathered}
\{x\in (-\infty,0]^2\mid H_2=0\}=\{G_2(t)\mid t\leq 0\}=
\{0\}\cup \{x_{QS}^\gl\mid \gl>0\}
\\(\text{resp., } \quad 
\{x\in (-\infty,0]^2\mid H_1=0\}=\{G_1(t)\mid t\leq 0\}=
\{0\}\cup \{x_{RS}^\gl\mid \gl>0\}). 
\end{gathered}
\]

\bproof Note first that $\psi(0)=\chi_i(0)=0$ for $i=1,2$. 

Let $\gl>0$. Note by \erf{gl to mu} that 
\[
x_{QS}^\gl=G_2\circ p_1^{-1}\circ g_1\circ f_1^{-1}(\gl) \ \ \AND \ \ 
x_{RS}^\gl=G_1\circ p_2^{-1}\circ g_2\circ f_2^{-1}(\gl). 
\]
From the above, we have
\[
m_\gam\cdot x_{QS}^\gl=q_1\circ p_1^{-1}\circ g_1\circ f_1^{-1}(\gl) 
\ \ \AND \ \ n_\gam\cdot x_{RS}^\gl=q_2\circ p_2^{-1}\circ g_2\circ f_2^{-1}(\gl).  
\]
Since \ $
x_{QS}^\gl, \,x_{RS}^\gl, \,m_\gam,\, n_\gam\in(-\infty,0)^2$, 
we have \ $
m_\gam\cdot x_{QS}^\gl>0\ \AND \ n_\gam\cdot x_{RS}^\gl>0$. 
Hence, by the definition of $\chi_i$, we find that 
\[
\chi_1(x_{QS}^\gl)=\gl \ \ \AND \ \ \chi_2(x_{RS}^\gl)=\gl. 
\]
Moreover, Lemmas \ref{I^gl} and  \ref{psi=psi_012} ensure that
$\,\chi_1=\psi_1=\psi\,$ on the curve $\,H_2=0\,$ in $S$  
and $\,\chi_2=\psi_2=\psi\,$ on the curve $\,H_1=0\,$ in $S$.

We fix any $x\in P\cup Q\cup R$ and show that $\chi_i(x)\leq \psi(x)$ for $i=1,2$. 
Once this is done, since $\chi_i$ are convex on $\R^2$, we conclude that 
$\chi_i\leq \psi$ on $\R^2$ for $i=1,2$ and finish the proof.

We may henceforth assume that $x\not=0$. 
Choose $\gl>0$ so that $\psi(x)=\gl$, and  
note that $x\in I_P^\gl \cup I_Q^\gl \cup I_R^\gl$.  

Consider the case $x\in I_P^\gl$. Since $x\in(0,\infty)^2$ and $m_\gam, \, n_\gam\in (-\infty,0)^2$,  we have $m_\gam\cdot x<0$, $n_\gam\cdot x<0$, and hence, 
$\chi_i(x)=0\leq \psi(x)$ for $i=1,2$.

Recall the monotonicity of $\chi_i$: for any 
$\xi,\eta\in\R^2$, if $m_\gam\cdot \xi\leq m_\gam\cdot \eta$, then 
$\chi_1(\xi)\leq \chi_1(\eta)$ and, if  $n_\gam\cdot \xi\leq n_\gam\cdot \eta$, then 
$\chi_2(\xi)\leq \chi_2(\eta)$. 

Now, assume that $x\in I_Q^\gl$. 
Note that $x=x_{QS}^\gl+\tau (x_{PQ}^\gl-x_{QS}^\gl)$ for some $\tau\in[0,\,1]$ 
and that $x_{PQ}^\gl-x_{QS}^\gl\in (0,\infty)^2$, and hence, $m_\gam\cdot(x_{PQ}^\gl-x_{QS}^\gl)<0$,\ $n_\gam\cdot (x_{PQ}^\gl-x_{QS}^\gl)<0$.  
Moreover, we have
\[\bald
m_\gam \cdot x\leq m_\gam\cdot x_{QS}^\gl \ \ \ \AND \ \ \ 
n_\gam\cdot x\leq n_\gam\cdot x_{QS}^\gl. 
\eald\]
Similarly, we find that if $x\in I_R^\gl$, then 
\[m_\gam\cdot x\leq m_\gam\cdot x_{RS}^\gl, \quad 
n_\gam\cdot x\leq n_\gam\cdot x_{RS}^\gl. 
\] 
Hence, if $x\in I_Q^\gl$, then 
\[
\chi_1(x)\leq \chi_1(x_{QS}^\gl)=\gl=\psi(x),
\]
and, if $x\in I_R^\gl$, then 
\[
\chi_2(x)\leq \chi_2(x_{RS}^\gl)=\gl=\psi(x). 
\]

By Lemma \ref{comp-slope}, there are points $y_1=G_2(r_1)$ and $y_2=G_1(r_2)$, 
with $r_i<0$ for $i=1,2$,  such that $m_\gam\cdot x_{QS}^\gl=m_\gam\cdot y_2$, 
$n_\gam\cdot x_{RS}^\gl=n_\gam\cdot y_1$. Moreover, the lemma says that if we write 
$x_{QS}^\gl=G_2(t)$ and $x_{RS}^\gl=G_1(s)$, then 
$r_1\leq t<0$ and $r_2\leq s<0$.  The last two inequalities imply that
\[
x_{QS}^\gl-y_1,\, x_{RS}^\gl-y_2\in[0,\infty)^2,
\]
and hence, 
\[
n_\gam\cdot (x_{QS}^\gl-y_1)\leq 0 \quad\AND\quad 
m_\gam\cdot(x_{RS}^\gl-y_2)\leq 0. 
\]
Consequently, if $x\in I_Q^\gl$, then $n_\gam\cdot x\leq n_\gam\cdot y_1$ and
\[\bald
\chi_2(x)& \leq \chi_2(y_1)= f_2\circ g_2^{-1}\circ p_2\circ q_2^{-1}(n_\gam \cdot y_1) 
= f_2\circ g_2^{-1}\circ p_2\circ q_2^{-1}(n_\gam \cdot x_{RS}^\gl) 
\\&=\chi_2(x_{RS}^\gl)=\gl=\psi(x),
\eald
\]
and  if $x\in I_R^\gl$, then $m_\gam\cdot x\leq m_\gam\cdot y_2$ and 
\[\bald
\chi_1(x)&\leq \chi_1(y_2) 
=f_1\circ g_1^{-1}\circ p_1\circ q_1^{-1}(m_\gam \cdot x_{QS}^\gl) 
= \chi_1(x_{QS}^\gl) =\gl=\psi(x).
\eald
\]
The proof is complete. 
\eproof

We remark that we may choose $\gam>0$ as small as we wish in 
\erf{cond-slope}. 

\begin{lemma} \label{normal}Let $\gl>0$ and $L_\gl:=\{x\in\R^2 \mid \psi(x)\leq \gl\}$. Let $y\in\pl L_\gl$.  We select $\gam>0$ sufficiently small. Then, 
\[ \bald
 & N(y,L_\gl)\subset \cone\{v_1, n_\gam\} \ \ & \IF \ y_1<h_1(y_2),
\\&
N(y,L_\gl)\subset \cone\{v_2, m_\gam\} \ \ &\IF \ y_2<h_2(y_1),
\\&
N(y,L_\gl)\subset \cone\{v_0,v_2\} \ \ &\IF \ y_1>h_1(y_2), 
\\& 
N(y,L_\gl)\subset \cone\{v_0,v_1\} \ \ &\IF \ y_2>h_2(y_1). 
\eald
\]
\end{lemma}

We formulate a simple lemma for the proof of the lemma above.

\begin{lemma}\label{cone} Let $v,w,p,q\in\R^2$. Assume that 
\[
v\cdot p=w\cdot q=0, \quad v\cdot q<0, \  \ \AND \ \ w\cdot p<0.
\]
Set \ $
C=\{x\in\R^2\mid v\cdot x\leq 0,\ w\cdot x\leq 0\}$ \ and \ $C^*=\{\xi\in\R^2 
\mid \xi\cdot x\leq 0 \ \FORALL x\in C\}$. 
Then, 
\[
C=\cone\{p,q\} \quad\AND\quad C^*=\cone\{v,w\}. 
\]
\end{lemma}

\bproof It is easily seen that the vectors $p$ and $q$ are linearly independent. 
For any $x=sp+tq$, with $s,t\in\R$, we have
\[
v\cdot x=t v\cdot q \ \ \AND \ \ w\cdot x=s w\cdot p,
\]
and, since $v\cdot q<0$ and $w\cdot p<0$, we find that 
$x\in C$ if and only if $s\geq 0$ and $t\geq 0$. Hence, we have
\[
C=\cone\{p,q\}.
\]

Next, note first that $v,w\in C^*$ and $C^*$ is a closed convex cone. Hence, 
we have $\cone\{v,w\}\subset C^*$. To show that $\cone\{v,w\}=C^*$, 
we suppose to the contrary that $C^*\stm \cone\{v,w\}\not=\emptyset$. Select 
$u\in C^*\stm\cone\{v,w\}$. By the Hahn-Banach theorem, there is $x_0\in\R^2$
such that $v\cdot x_0\leq 0$, $w\cdot x_0\leq 0$, and $u\cdot x_0>0$, which 
can be stated that $x_0\in C$ and $u\cdot x_0>0$. These contradict 
that $u\in C^*$. 
\eproof

\bproof[Proof of Lemma \ref{normal}]   We consider only 
the case when $y_1<h_1(y_2)$. The other cases are treated similarly, whose 
details are left to the reader.

Set \ 
$
C=\{x\in\R^2\mid v_1\cdot x\leq 0,\ n_\gam\cdot x\leq 0\},
$ 
which is a closed convex cone with vertex at the origin.  
We prove that 
\beq\label{inc-NC*}
N(y,L_\gl)\subset C^*,
\eeq
where \ $
C^*=\{\xi\in\R^2\mid \xi\cdot x\leq 0 \ \FORALL x\in C\}.$ 
Once this is done, thanks to 
Lemma \ref{cone}, we conclude the proof. 

To check \erf{inc-NC*}, we set $p=x_{PQ}^\gl$, $q=x_{RS}^\gl$, 
and $D=\{x\in\R^2\mid v_1\cdot(x-p)\leq 0,\ n_\gam\cdot (x-q)\leq 0\}$, 
and note by the definition of $L_\gl$ that \ $
L_\gl\subset D$ \ and \ $p,q\in\pl L_\gl.$

Since $v_1$ and $n_\gam$ are linearly independent, there is a unique solution 
$x\in\R^2$ of the system \ $v_1\cdot(x-p)=n_\gam\cdot (x-q)=0$. 
Let $r$ denote the unique solution $x$ of the above system. Note that 
\[\bald 
D&=\{x\in\R^2\mid v_1\cdot(x-r)\leq 0,\ n_\gam\cdot (x-r)\leq 0\}
=r+C.
\eald
\]

The proof of Lemma 
\ref{4 func} easily extends to show that 
the functions 
\[
t\mapsto v_1\cdot G_1(t) \ \ \AND \ \ t\mapsto n_\gam\cdot G_1(t)
\]
are, respectively, uniformly increasing and decreasing on $\R$. 
The characterization of the points $p,q$ due to Proposition \ref{four-points}
shows that $v_1\cdot q<0<v_1\cdot p$ and $n_\gam\cdot p <0<n_\gam\cdot q$. 
Also, since $v_1\cdot e_1<0$,  
$n_\gam\cdot e_1<0$,  and $y_1<h_1(y_2)$ by assumption, we have 
\[\bald
v_1\cdot p\geq v_1\cdot y>v_1\cdot G_1(y_2) \ \ \AND \ \ 
n_\gam\cdot q\geq n_\gam\cdot y>n_\gam\cdot G_1(y_2), 
\eald \]
which assures that $q_2<y_2<p_2$.

Noting that 
\beq\label{vnpq}
v_1\cdot (p-r)=0>v_1\cdot (q-r) \ \ \AND \ \ n_\gam\cdot(q-r)=0>n_\gam\cdot(p-r),
\eeq
we invoke Lemma \ref{cone}, to find that 
\beq\label{cones}
C=\cone\{p-r,q-r\} \ \ \AND \ \ C^*=\cone\{v_1,n_\gam\}. 
\eeq

Since $y-r\in C$, there are unique $s,t\in[0,\infty)$ such that 
$y-r=s(p-r)+t(q-r)$. 
We claim that $s+t<1$. Indeed,  
if we write $e_1=\gs (p-r)+\tau(q-r)$, 
with $\gs,\tau\in\R$, we find by \erf{vnpq} that $\gs>0$ and $\tau>0$.  
Since $q_2<y_2<p_2$, we have $G_1(y_2)\in L_\gl$. 
Note that $G_1(y_2)=y+(h_1(y_2)-y_1)e_1$. There exists $\gth>0$ 
such that $y+\gth e_1\in [q,p]$, which reads that the point
\[
y-r+\gth e_1=(s+\gth \gs)(p-r)+(t+\gth\tau)(q-r)
\]
lies in the interval $[q-r,\,p-r]$. This yields 
\[
s+\gth\gs+ t+\gth\tau=1,
\]
which ensures that $s+t<1$. 

Now, let $\xi\in N(y,L_\gl)$, which implies that 
\[
\xi\cdot(q-y)\leq 0 \ \ \AND \ \ \xi\cdot (p-y)\leq 0. 
\]
Since $y-r=s(p-r)+t(q-r)$, we have
\[\bald
& 0\leq \xi\cdot(y-p)=\xi\cdot(r-p)+s\xi\cdot(p-r)+t\xi\cdot(q-r), 
\\& 0\leq \xi\cdot(y-q)=\xi\cdot(r-q)+s\xi\cdot(p-r)+t\xi\cdot(q-r).
\eald
\]
Combining these with the inequalities $s+t<1,\,s\geq 0,\,t\geq 0$, we infer that 
\[
\xi\cdot (p-r)\leq 0 \ \ \AND \ \ \xi\cdot(q-r)\leq 0,
\]
which imply together with the first identity of \erf{cones} that $\xi\in C^*$.  
The proof is complete.
\eproof 

\begin{lemma} \label{psi^2} For any $x\in\R^2$, 
 \[
\pl \psi^2(x)=2\psi(x)\pl \psi(x). 
\]
\end{lemma}

\bproof  If $p\in \pl\psi(x)$, then, as $\R^2\ni y\to 0$, 
\[\bald
\psi^2(x+y)-\psi^2(x)&=(\psi(x+y)+\psi(x))(\psi(x+y)-\psi(x))
\\&\geq (\psi(x+y)+\psi(x)) p\cdot y 
=2\psi(x)p\cdot y+O(|y|^2),
\eald\]
which assures that $2\psi(x)p\in \pl \psi^2(x)$. If $p\in \pl\psi^2(x)$ and $x\not=0$, then, as $y\to 0$, 
\[\bald
(\psi(x+y)+\psi(x))(\psi(x+y)-\psi(x))=\psi^2(x+y)-\psi^2(x)
\geq p\cdot y,
\eald\]
and 
\[\bald
\psi(x+y)-\psi(x)
\geq \fr{p\cdot y}{\psi(x+y)+\psi(x)}
=\fr{p}{2\psi(x)}\cdot y+O(|y|^2),
\eald\]
which reads that $p/(2\psi(x))\in \pl \psi(x)$. On the other hand, we clearly  have 
$\pl \psi^2(0)=\{0\}=2\psi(0)\pl \psi(0)$. 
\eproof 

\begin{lemma}\label{bounds-Dpsi} For some constants $0<\gd\leq C<\infty$, 
\[
\gd\leq |p|\leq C \ \ \FORALL p\in \pl \psi(x) \ \AND \ x\in\R^2\stm\{0\}.
\] 
\end{lemma}

\bproof According to Lemma \ref{psi l-growth}, there exist constants 
$0<\gd_0\leq C_0<\infty$ such that 
\[
\gd_0|x|\leq \psi(x)\leq C_0|x| \ \ \FORALL \ x\in\R^2.
\]
Let $x\in\R^2\stm\{0\}$ and $p\in\pl \psi(x)$.  For any $y\in\R^2$, 
we have 
\beq\label{bounds-Dpsi.1}
\psi(x+y)-\psi(x)\geq p\cdot y. 
\eeq
Put $y=-x$, to obtain
\[
\gd_0|x|\leq \psi(x)\leq p\cdot x\leq |p||x|, 
\]
and hence, $|p|\geq\gd_0$. Next, set $y=|x|p/|p|$, to find that
\[
|p||x|\leq \psi\left(x+\tfrac{|x|p}{|p|}\right)\leq 2C_0|x|,
\] 
which yields that $|p|\leq 2C_0$. Thus, we have \ $\gd_0\leq|p|\leq 2C_0$. 
\eproof

\section{$C^{1,1}$ test function}

In the previous section, we have constructed a convex function $\psi$ on $\R^2$ such that for some constants $0<\gd\leq C<\infty$, 
\beq \label{deltaC1}
\gd|x|\leq \psi(x)\leq C|x| \ \ \FORALL \ x\in \R^2,
\eeq
and 
\beq\label{deltaC2}
\gd\leq |p|\leq C \ \ \FORALL \ p\in \pl\psi(x),\, x\in\R^2\stm\{0\}.
\eeq 
(See Lemmas \ref{psi l-growth} and \ref{bounds-Dpsi} for the estimates above.)
More importantly, thanks to Lemma \ref{normal}, if we set for $\gl>0$,
\[
L_\gl:=\{p\in\R^2\mid \psi(p)\leq\gl\}, 
\]
then for any $q\in \pl L_\gl$, we have 
\beq \label{normal2} N(q,L_\gl)\subset \bcases
\cone\{v_1,n_\gam\}  \ & \IF\ q_1<h_1(q_2),\\
\cone\{v_2,m_\gam\} & \IF\ q_2<h_2(q_1),\\
\cone\{v_0,v_2\} & \IF\ q_1>h_1(q_2), \\
\cone\{v_0,v_1\} & \IF\ q_2>h_2(q_1).
\ecases
\eeq

For $\ep>0$, we set 
\[
\eta_\ep(p)=\psi^2(p)+\ep|p|^2,
\]
and
\beq \label{def-phi}
\phi_\ep(x)=\max_{p\in\R^2}(x\cdot p-\eta_\ep(p)).
\eeq
Note that $\phi_\ep$ is convex on $\R^2$,
\[\phi_\ep\in C^{1,1}(\R^2), \ \ \AND \ \ 
\gd'|x|^2\leq \phi_\ep(x)\leq C'|x|^2 \ \ \FORALL \ x\in\R^2,
\]
where $0<\gd'\leq C'<\infty$ are constants.
Indeed, by the definition of $\phi_\ep$, it is clear that $\phi_\ep$ is convex on $\R^2$. Writing 
\beq \label{inf-conv}
\phi_\ep(x)=\tfrac{1}{4\ep}|x|^2-\min_{p\in\R^2}\left(\psi^2(p)
+\ep\left|p-\tfrac{1}{2\ep}x\right|^2\right),
\eeq
we find that $\phi_\ep$ is semi-concave on $\R^2$. Hence, $\phi_\ep\in C^{1,1}(\R^2)$.  (Notice that the function $x\mapsto -\phi_\ep(x)+|x|^2/(4\ep)$ is the inf-convolution of the function $\psi^2$.) Finally, note that, since 
\[
(\ep+\gd^2)|p|^2\leq \eta_\ep(p)\leq (\ep+C^2)|p|^2, 
\]
we have 
\[
\phi_\ep(x)
\bcases
\leq \max_{p\in\R^2}\left(x\cdot p-(\ep+C^2)|p|^2\right)= \tfrac{1}{4(\ep+C^2)}|x|^2& \\
\geq \max_{p\in\R^2}\left(x\cdot p-(\ep+\gd^2)|x|^2\right)=\tfrac{1}{4(\ep+\gd^2)|x|^2},& 
\ecases
\]
and conclude that for all $x\in\R^2$,  
\[
\tfrac{1}{4(\ep+C^2)}|x|^2\leq \phi_\ep(x)\leq \tfrac{1}{4(\ep+\gd)}|x|^2.
\]

\begin{theorem} \label{exist-phi} Let $\ep>0$ be sufficiently small. 
The function $\phi_\ep$ defined by \erf{def-phi} satisfies: 
\beq \label{cond-phi_ep}
H_1(D\phi_\ep(x)) \bcases
\geq 0 & \IF \ x_1\leq 0,\\
\leq 0 &\IF\ x_1\geq 0,
\ecases  \ \ \ \AND \ \ \ 
H_2(D\phi_\ep(x))\bcases
\geq 0 & \IF \ x_2\leq 0,\\
\leq 0 & \IF \ x_2\geq 0.
\ecases
\eeq
\end{theorem}

A possible choice of the constant $\ep>0$ in the theorem above is 
\beq \label{cond-ep}
0<\ep<\fr{\ep_0 \gd^2}{\max\{|v_0|,|v_1|,|v_2|,|m_\gam|,|n_\gam|\}}, 
\eeq 
where $\gd$ is the constant appearing in \erf{deltaC1} and \erf{deltaC2} 
and 
\beq\label{ep_0}
\ep_0=\min\{v_0\cdot e_1,v_0\cdot e_2, 
-v_1\cdot e_1, v_1\cdot e_2, v_2\cdot e_1, -v_2\cdot e_2, -m_\gam\cdot e_2,-n_\gam \cdot e_1\}. 
\eeq
Since $m_\gam\cdot e_1=-\gam=n_\gam\cdot e_2<0$, we find by \erf{cond-v_i} that   
$\ep_0>0$. 

\bproof We assume henceforth that $\ep$ satisfies \erf{cond-ep}, and prove that 
\erf{cond-phi_ep} holds. 
As before, we set $ L_\gl=\{q\in\R^2\mid \psi(q)\leq\gl\}$ for $\gl>0$. 
Fix any $x\in\R^2$. Assume first $x_1\leq 0$ and show that
$H_1(D\phi_\ep(x))\geq 0$. We argue by contradiction. Set $p=D\phi_\ep(x)$ and suppose that
$
H_1(p)=-p_1+h(p_2)<0.
$

We have $p_1>h_1(p_2)$ and $p\not=0$. 
By the convex duality, using  Lemma \ref{psi^2}, we have 
\[
x\in \pl \eta_\ep(p)=2\psi(p)\pl\psi(p)+2\ep p, 
\]
that is,  $x-2\ep p\in 2\psi(p)\pl\psi(p)$. 
Set 
$\gl=\psi(p), $
and note that $p\in \pl L_\gl$, $x-2\ep p\in N(p,L_\gl)$, and 
$|x-2\ep p|\geq 2\gl\gd\geq 2\gd^2|p|$ by \erf{deltaC1} and \erf{deltaC2}. 
Since $p_1>h_1(p_2)$, we have, by \erf{normal2}, 
\[
N(p,L_\gl)\subset \cone\{v_0,v_2\}. 
\] 
Accordingly, there are $t,s\geq 0$, with $t+s>0$, such that
$x-2\ep p=tv_0+sv_2$.  We have 
\[
2\gd^2|p|\leq |x-2\ep p|\leq \max\{|v_0|,|v_2|\}(t+s),
\]
and also,  by \erf{ep_0}, 
\[
x_1\geq e_1\cdot (x-2\ep p)-2\ep |p| \geq \ep_0(t+s)-2\ep|p|. 
\]
Hence, 
\[
\max\{|v_0|,|v_2|\}x_1\geq \max\{|v_0|,|v_2|\}(\ep_0(t+s)-2\ep |p|)
 \geq 2\left(\ep_0\gd^2-\ep \max\{|v_0|,|v_2|\} \right)|p|.
\]
This, together with \erf{cond-ep}, gives a contradiction that $x_1>0$. 

Next, we assume that $x_1\geq 0$ and show that  $H_1(D\phi_\ep(x))\leq 0$. 
The proof in this case is similar to the above, and we argue by contradiction. 
Suppose that $H_1(p)> 0$, that is, $p_1<h_1(p_2)$, where $p=D\phi_\ep(x)$.    
We have $x-2\ep p\in 2\gl \pl\psi(p)$, where $\gl=\psi(p)$, and by \erf{normal2},
\[
\pl\psi(p)\subset N(p,L_\gl)\subset \cone\{v_1,n_\gam\}. 
\]
We choose $t,s\geq 0$, with $t+s>0$, so that $x-2\ep p=t v_1+s n_\gam$ and note that
\[
2\gd^2|p|\leq |x-2\ep p|\leq \max\{|v_1|,|n_\gam|\}(t+s),
\] 
and, by \erf{ep_0},
\[
x_1\leq e_1\cdot(x-2\ep p)+2\ep|p|
\leq -\ep_0(t+s)+2\ep|p|.
\]
Consequently, we have
\[
 \max\{|v_1|,|n_\gam|\} x_1\leq -2\left(\ep_0\gd^2 -\ep\max\{|v_1|,|n_\gam|\} \right)|p|,
\]
which yields a contradiction. 

Now, we assume that $x_2\leq 0$. Similarly to the above, 
we set $p=D\phi_\ep(x)$ and  suppose to the contrary that $H_2(p)< 0$, that is, $p_2>h_2(p_1)$.    
We have $x-2\ep p\in 2\gl \pl\psi(p)$, with $\gl=\psi(p)$, and by \erf{normal2},
\[
x-2\ep p\in N(p,L_\gl)\subset \cone\{v_0,v_1\}. 
\]
We choose $t,s\geq 0$, with $t+s>0$, so that $x-2\ep p=t v_0+s v_1$ and compute as before that
\[
2\gd^2|p|\leq |x-2\ep p|\leq \max\{|v_0|,|v_1|\}(t+s),
\] 
and
\[
x_2\geq e_2\cdot(x-2\ep p)-2\ep|p|
\geq \ep_0(t+s)-2\ep|p|.
\]
Combining these two, we have
\[
 \max\{|v_0|,|v_1|\} x_2\geq 2\left(\ep_0\gd^2 -\ep\max\{|v_0|,|v_1|\} \right)|p|.
\]
This yields a contradiction, which shows that $H_2(p)\geq 0$.

Finally, we assume that $x_2\geq 0$ and show that  $H_2(p)\leq 0$, where $p=D\phi_\ep(x)$.  As always, we argue by contradiction and, hence, suppose that $H_2(p)> 0$, that is, $p_2<h_2(p_1)$.     
We have $x-2\ep p\in 2\gl \pl\psi(p)$, where $\gl=\psi(p)$, and by \erf{normal2},
\[
x-2\ep p\in  N(p,L_\gl)\subset \cone\{v_2,m_\gam\}. 
\]
We choose $t,s\geq 0$, with $t+s>0$, so that $x-2\ep p=t v_2+s m_\gam$ and note that
\[
2\gd^2|p|\leq |x-2\ep p|\leq \max\{|v_2|,|m_\gam|\}(t+s),
\] 
and
\[
x_2\leq e_2\cdot(x-2\ep p)+2\ep|p|
\leq -\ep_0(t+s)+2\ep|p|.
\]
Hence, we have
\[
 \max\{|v_2|,|m_\gam|\} x_2\leq -2\left(\ep_0\gd^2 -\ep\max\{|v_2|,|m_\gam|\} \right)|p|,
\]
and get a contradiction.  
\eproof 

\def\th{\tilde h}

\section{Comparison principle} \label{comparison}

As an application of Theorem \ref{exist-phi}, 
we establish a comparison theorem for \erf{N-prob} in this section. 

Let $\gO:=(0,\infty)^2$, as in the previous sections, and we reformulate Theorem 
\ref{exist-phi} by removing the normalization \erf{zero-norm}. 

Under the assumption \erf{cond-H_i}, the functions $h_i$ are both Lipschitz continuous on $\R$, with Lipschitz bound $\ga_0:=\max\{\ga_1,\ga_2\}\in[0,1)$, and, if we define the function $h\mid \R^2\to\R^2$ by $h(x)=(h_1(x_2),h_2(x_1))$, then 
$h$ is Lipschitz continuous, with Lipschitz bound $\ga_0$. By the standard fixed point theorem, there exists a unique fixed point $a=(a_1,a_2)\in\R^2$ of the mapping $h$. 
In particular, we have 
\beq \label{fixed point} h_1(a_2)=a_1\quad\AND\quad h_2(a_1)=a_2. 
\eeq
We define the functions 
$\tH_i\mid \R^2\to\R$, $\tilde h_i\mid \R\to\R$, $i=1,2$, by
\beq\label{tilde H}
\tH_i(p)=H_i(p+a),\quad \th_1(t)=-a_1+h_1(t+a_2), \quad \th_2(t)=-a_2+h_2(t+a_1),
\eeq   
and note that for all $p\in\R^2$ and $i=1,2$, 
\[
\tH_1(p)=-p_1+\th_1(p_2), \quad \tH_2(p)=-p_2+\th_2(p_1), 
\quad \th_i(0)=0. 
\]
It is easily checked that \erf{cond-H_i} holds with $\tH_i$ replacing $H_i$.

\begin{corollary} \label{exist-f}  Assume \erf{cond-H_i}. Let $a\in\R^2$ be the point such that \erf{fixed point} holds. Define $\tH_i$ and $\th_i$, $i=1,2$, 
by \erf{tilde H}, and assume that \erf{cond-slope} is satisfied for $\tH_i$ 
in place of $H_i$, $i=1,2$, and a choice of vectors $v_i$, $i=0,1,2$. 
Then, there exists a convex function $f\in C^{1,1}(\R^2)$ such that 
\beq \label{f1}
H_1(Df(x))
\bcases 
\geq 0\ \ &\IF\ x_1\leq 0, \\
\leq 0 &\IF\ x_1\geq 0, 
\ecases
\quad\AND \quad
H_2(Df(x))
\bcases
\geq 0 \ \ &\IF\ x_2\leq 0, \\
\leq 0 &\IF\ x_2\geq 0,
\ecases
\eeq
and for some positive constants $\gd, C$, 
\beq \label{f2}
\gd|x|^2-C\leq f(x)\leq C(|x|^2+1) \ \ \ \FORALL x\in\R^2. 
\eeq
\end{corollary}

\bproof We apply Theorem \ref{exist-phi}, to choose a function $\phi\in C^{1,1}(\R^2)$ which satisfies the conditions \erf{f1}, with $H_i$ replaced by $\tH_i$, and 
\erf{f2} for some constants $0<\gd<C<\infty$.  
We set $f(x)=\phi(s)+a\cdot x$ for $x\in\R^2$, observe that for all $x\in\R^2$ and $i=1,2$,
\[
H_i(Df(x))=H_i(D\phi(x)+a)=\tH_i(D\phi(x)),
\]
and conclude that $f$ satisfies \erf{f1}. It is clear that $f$ satisfies \erf{f2}. 
\eproof

As an application of Corollary \ref{exist-f}, we present a comparison theorem for bounded sub and supersolutions 
of \erf{N-prob}. 

We write $\bS^n$ for the space of real $n\tim n$ symmetric matrices and say that 
a function $\go: [0,\infty)\to[0,\infty)$ is a modulus if it is continuous 
and nondecreasing on $[0,\infty)$. We denote by $I_n$ the unit matrix in $\bS^n$. 
For $X,Y\in\bS^n$, the inequality $X\leq Y$ means that $Y-X$ is positive  semidefinite.  The choice of norm for $X\in\bS^n$ here is 
$|X|=\max\{ |X\xi|\mid \xi\in\R^n, \,|\xi|=1\}$. 

Let $F$ be a continuous, real-valued function on $\bS^2\tim\R^2\tim\R\tim\ol\gO$.
Assume in the spirit of \cite{CIL1992} that $F$ satisfies the following three conditions
\erf{F1}--\erf{F3}.  
\renewcommand{\theenumi}{F\arabic{enumi}}
\renewcommand{\labelenumi}{(\theenumi)}

\begin{enumerate}
\item
\label{F1}   
There is a constant $\mu>0$ such that for each $(X,p,x)\in\bS^2\tim\R^2\tim\gO$, the function \ $u\to F(X,p,u,x)-\mu u$\ is nondecreasing on $\R$.
\item 
\label{F2}  
There is a constant $L>0$ such that for every $(X,p),(Y,q)\in \bS^2\tim \R^2$ and 
$(u,x)\in\R\tim\gO$,  
\[  
|F(X,p,u,x)-F(Y,q,u,x)|\leq L(|X-Y|+|p-q|). 
\] 
\item \label{F3}  
There is a modulus $\go$ such that for any $X,Y\in \bS^2, x,y\in\gO$, $(p,u)\in\R^2\tim\R$, 
and $\ga>0$, if 
\[
-\ga \bmat I_2 & 0\\ 0&I_2\emat \leq \bmat X&0\\0&-Y\emat 
\leq \ga\bmat I_2&-I_2\\ 
-I_2&I_2\emat,
\]
then 
\[
F(X,p,u,x)-F(Y,p,u,y)\geq -\go(\ga|x-y|^2+|x-y|(|p|+1)). 
\]
\end{enumerate}  

Note (see \cite{CIL1992}) that if $F\in C(\bS^2\tim\R^2\tim\R\tim\gO)$ satisfies \erf{F3}, then $F$ is degenerate elliptic. The degenerate ellipticity means here that 
if $X,Y\in\bS^2$ and $X\leq Y$, then $F(X,p,u,x)\geq F(Y,p,u,x)$ for all $(p,u,x)\in
\R^2\tim\R\tim\gO$. 

We remark that a large class of elliptic operators arising as the dynamic programming 
equations (Hamilton-Jacobi-Bellman or Hamilton-Jacobi-Isaacs equations) in stochastic optimal control or differential games satisfy \erf{F1}--\erf{F3}.

\begin{theorem}\label{comp.thm} Let $F\in C(\bS^2\tim\R^2\tim\R\tim\gO,\R)$ satisfy \erf{F1}--\erf{F3}. Let $H_i\in C(\R^2,\R)$, $i=1,2$, satisfy \erf{cond-H_i}. 
Let $a\in\R^2$ satisfy \erf{fixed point} and define $\tH_i$, $i=1,2$, by \erf{tilde H}. 
Assume that \erf{cond-slope} holds for $\tH_i$, $i=1,2$, and some $v_i\in\R^2$, $i=0,1,2$, satisfying \erf{cond-v_i}--\erf{cond-v_2}. Let $v,w\mid\ol\gO\to \R$
be, respectively, a lower semicontinuous viscosity subsolution and 
upper semicontinuous viscosity supersolution of \erf{N-prob}. Moreover, assume that 
$v,w$ are bounded on $\ol\gO$. Then, $v\leq w$ on $\ol\gO$.  
\end{theorem}

We state an existence theorem of a (continuous) viscosity solution 
of \erf{N-prob}, where the boundedness of the function $x\mapsto F(0,0,0,x)$ 
assures the existence of a bounded solution.

\begin{corollary}\label{exist} In addition to the hypotheses of Theorem \ref{comp.thm}, 
assume that the function $\gO\ni x\mapsto F(0,0,0,x)$ is bounded. Then, there exists a unique viscosity solution $u\in C(\ol\gO)$ of \erf{N-prob}. 
\end{corollary}

\bproof The uniqueness assertion is an obvious consequence of 
Theorem \ref{comp.thm}. 

Choose a nonnegative function $\eta\in C^2([0,\infty))$ such that 
\[
\eta'(0)=-1 \quad \AND\quad\max\{\eta(t), |\eta'(t)|,|\eta''(t)|\}\leq 1 \ \ \FORALL t\geq 0,
\]
set $\psi(x)=\eta(x_1)+\eta(x_2)$ for $x\in\ol\gO$, and observe that for any 
$x\in\ol\gO$ and constant $C>0$, if $x_1=0$, 
\[
H_1(CD\psi(x))\geq -C\eta'(0)+h_1(C\eta'(x_2)\leq h_1(0)+C(1-\ga_1),
\]
and
\[
H_1(-CD\psi(x))\leq C\eta'(0)+h_1(-C\eta'(x_2))\leq h_1(0)-C(1-\ga_1). 
\]
If $x_2=0$, then
\[
H_2(CD\psi(x))\geq -C\eta'(0)+h_2(C\eta'(x_1)\geq h_2(0)+C(1-\ga_2),
\]
and
\[
H_2(-CD\psi(x))\leq C\eta'(0)+h_2(-C\eta'(x_2))\leq h_2(0)-C(1-\ga_2). 
\]
For instance, we may choose the function $\eta(t)=e^{-t}$. 

Choose a constant $C_0>0$ so that $|F(0,0,0,x)|\leq C_0$ for all $x\in\gO$. 
Let $C_1,C_2$ be positive constants to be fixed later, and set 
\[
v(x)=-C_1\psi(x)-C_2 \ \ \ \AND \ \ \ w(x)=C_1\psi(x)+C_2\quad\FOR x\in\ol\gO.
\]
Compute that for any $x\in\ol\gO$,
\[\bald
F(D^2v(x),Dv(x),v(x),x)
&\leq F(0,0,v(x),x)-\mu v(x)+\mu v(x)
\\&\quad +LC_1(|\diag\{\eta''(x_1),\eta''(x_2)\}|+2)
\\&\leq C_0-\mu C_2+3LC_1, 
\eald\]
if $x_1=0$, then
\[
H_1(Dv(x))=H_1(-C_1D\psi(x))\leq h_1(0)-C_1(1-\ga_1), 
\]
and if $x_2=0$, then
\[
H_2(Dv(x))\leq h_2(0)-C_1(1-\ga_2).
\]
Similarly, we have for any $x\in\ol\gO$,
\[
F(D^2w(x),Dw(x),w(x),x)
\geq F(0,0,0,x)+\mu w(x)-3LC_1
\geq -C_0+\mu C_2-3LC_1, 
\]
if $x_1=0$, then
\[
H_1(Dw(x))\geq h_1(0)+C_1(1-\ga_1), 
\]
and if $x_2=0$, then
\[
H_2(Dw(x))\leq h_2(0)+C_1(1-\ga_2).
\]
We fix $C_1$ so that 
\[
|h_i(0)|\leq C_1(1-\ga_i) \ \ \ \FOR i=1,2,
\]
and then $C_2$ so that  
\[
\mu C_2\geq C_0+3LC_1,
\]
to find that $v$ and $w$ are, respectively, a classical (and hence, viscosity) 
subsolution and supersolution of \erf{N-prob}. 
If we set
\[
u(x)=\sup\{\phi(x)\mid v\leq \phi\leq w \ON \ol\gO,\ \phi \text{ viscosity subsolution of } \erf{N-prob}\} \ \ \FOR x\in\ol\gO,
\]
then, according to the Perron method (see \cite{CIL1992}), 
$u$ is an upper semicontinuous subsolution of \erf{N-prob} 
and the lower semicontinuous envelope $u_*$ of $u$ is a supersolution of \erf{N-prob}. 
Theorem \ref{comp.thm} assures that $u\leq u_*$ on $\ol\gO$, which implies that 
$u\in C(\ol\gO)$. Since $v\leq u\leq w$ on $\ol\gO$, it is clear that $u$ is bounded on $\ol\gO$.  
\eproof

\bproof[Proof of Theorem \ref{comp.thm}] We fix a nonnegative function $\eta\in C^2([0,\infty))$ so that 
\[
\eta'(0)=-1, \quad \max\{|\eta'(t)|,|\eta''(t)|\}\leq 1 \ \ \FORALL t\geq 0, \ \ \AND \ \ 
\lim_{t\to\infty} \eta(t)=+\infty.
\]
Let $\ep>0$, $C_1>0$, and set for $x\in\ol\gO$.
\[
\psi(x)=\eta(x_1)+\eta(x_2), \quad 
\tilde v(x)=v(x)-\ep(\psi(x)+C_1) \ \ \AND \ \ \tilde w(x)=w(x)+\ep(\psi(x)+C_1).
\]

We first observe that if $C_1$ is sufficiently large, then, for any $\ep>0$, $\tilde v$ and $\tilde w$ are a ``strict'' subsolution and supersolution of \erf{N-prob}, respectively. 
To check this, we compute slightly in an informal manner that for $x\in\gO$, 
\[
\bald
|D^2\psi(x)|&= |\diag\{\eta''(x_1),\eta''(x_2)\}|\leq 2,
\quad |D\eta(x)|\leq 1, 
\\F(D^2\tilde v,D\tilde v,\tilde v,x)
&\leq F(D^2 v,Dv,\tilde v,x)-\mu \tilde v+\mu \tilde v+3\ep\,L
\\& \leq F(D^2 v,Dv,v,x)+3\ep\,L-\mu \ep C_1\leq -\ep(\mu C_1-3L), 
\eald
\]
for $x\in \pl_1\gO$, 
\[
H_1(D\tilde v(x))=H_1(Dv(x)-\ep D\psi(x))\leq -\ep(1-\ga_1),
\]
and for $x\in\pl_2\gO$,
\[
H_2(D\tilde v(x))=H_2(Dv(x)-\ep D\psi(x))\leq -\ep(1-\ga_1).
\]
Similarly, if $x\in\gO$, 
\[
\bald
F(D^2\tilde w,D\tilde w,\tilde w,x)
\geq F(D^2 w,Dw,w,x)-3\ep\,L+\mu \ep C_1\geq \ep(\mu C_1-3L), 
\eald
\]
if $x\in \pl_1\gO$, 
\[
H_1(D\tilde w(x))=H_1(Dw(x)+\ep D\psi(x))\geq \ep (1-\ga_2),
\]
and if $x\in\pl_2\gO$,
\[
H_2(D\tilde w(x))=H_2(Dw(x)+\ep D\psi(x))\geq \ep (1-\ga_2). 
\]

Fixing  $C_1$ by $\mu C_1=3L+1$ and setting $\rho:=\min_{i=1,2}(1-\ga_i)$, 
we easily find that $\tilde v$ and $\tilde w$ are respectively a subsolution of
\beq\label{strict-sub}
\left\{\bald
&F(D^2u,Du,u,x)=-\ep\ \ \IN \gO,
\\&H_1(Du)=-\ep\rho \ \ \ON \pl_1\gO,\quad H_2(Du)=-\ep\rho \ \ \ON \pl_2\gO,
\eald \right. 
\eeq
and supersolution of
\beq\label{strict-super}
\left\{\bald
&F(D^2u,Du,u,x)=\ep \ \ \IN \gO,
\\&H_1(Du)=\ep\rho \ \ \ON \pl_1\gO,\quad H_2(Du)=\ep\rho \ \ \ON \pl_2\gO,
\eald \right. 
\eeq
Note also that 
\beq\label{at-infty}
\lim_{|x|\to \infty}\tilde v(x)=-\infty \ \ \AND \ \ 
\lim_{|x|\to \infty}\tilde w(x)=+\infty.
\eeq

It is enough to prove that $\max_{\ol\gO}(\tilde v-\tilde w)\leq 0$
 for all $\ep>0$. 

To do this, we argue by contradiction, and suppose that 
for some $\ep>0$,
\beq\label{positive}
\max_{\ol\gO}(\tilde v-\tilde w)>0.
\eeq
We fix $\ep>0$ so that \erf{positive} holds.  
Let $f$ be the function from Corollary \ref{exist-f}. 
Let $\ga>0$ and consider the function
\[
\tilde v(x)-\tilde w(y)-\ga^{-1} f(\ga(x-y)) 
\]
on $\ol\gO\tim\ol\gO$. In view of \erf{at-infty} and \erf{f2}, this function takes a maximum 
at some $(x_\ga,y_\ga)\in{\ol\gO\,}^2$. Moreover, the set 
$\{(x_\ga,y_\ga)\mid \ga>0\}$ is bounded. 

Since $f$ is a $C^{1,1}$ function, there is a constant $B>0$ such that $\xi\mapsto 
f(\xi)-B|\xi|^2/2$ is concave on $\R^2$, and hence,
\beq\label{semi-cv}
f(x)-f(\xi)\leq Df(\xi)\cdot (x-\xi)+\tfrac{B}{2}|x-\xi|^2 \ \ \FORALL x,\xi\in\R^2. 
\eeq
Consequently,  we have for any $x,y\in{\ol\gO\,}^2$, 
\[\bald
\tilde v(x)-\tilde w(y) &
-\tilde v(x_\ga)-\tilde w(y_\ga) 
\leq \ga^{-1}(f(\ga(x-y))-f(\ga(x_\ga-y_\ga)) 
\\& \leq Df(\ga(x_\ga-y_\ga))\cdot (x-y-((x_\ga-y_\ga)) 
+\tfrac{B \ga}{2}|x-y-(x_\ga-y_\ga)|^2,
\eald\] 
which implies that if we write $\,p_\ga=Df(x_\ga-y_\ga)\,$ and
\[
\psi_\ga(x,y)=p_\ga \cdot (x-y-(x_\ga-y_\ga)) 
+\tfrac{B \ga}{2}|x-y-(x_\ga-y_\ga)|^2 \ \ \FOR x,y\in\R^2, 
\]
then the function $(x,y)\mapsto \tilde v(x)-\tilde w(y) -\psi_\ga(x,y)$ has a maximum 
at $(x_\ga,y_\ga)$ on ${\ol\gO\,}^2$. Noting that 
\[
D^2\psi_\ga(x,y)=B\ga \bmat I_2&-I_2 \\ -I_2 & I_2\emat,
\] 
we deduce thanks to \cite[Theorem 3.2]{CIL1992} that
there are $X_\ga, Y_\ga\in\bS^2$ such that 
\begin{gather}\label{semi-max1}
(p_\ga, X_\ga)\in\ol{J}^{2,+}\tilde v(x_\ga),\quad
(p_\ga,Y_\ga)\in \ol{J}^{2,-}\tilde w(y_\ga), 
\\ \label{semi-max2}
-3B \ga\bmat I_2 &0 \\ 0& I_2\emat
\leq \bmat X_\ga & 0 \\ 0 & -Y_\ga\emat
\leq 3B \ga\bmat I_2&-I_2 \\ -I_2& I_2\emat,  
\end{gather}
where $\ol{J}^{2,+}$ and $\ol {J}^{2,-}$ 
are the ``closures'' of the second-order superjet and subjet, respectively (see \cite{CIL1992} for the definition of $\ol{J}^{2,\pm}$). 
It follows from \erf{semi-max1}, \erf{strict-sub}, and \erf{strict-super} that
\begin{gather}\label{st-visco1}
F(X_\ga, p_\ga, \tilde v(x_\ga),x_\ga)\leq -\ep \ \ \IF x_\ga\in\gO, \\ 
\label{st-visco2}
F(Y_\ga, p_\ga, \tilde w(y_\ga),y_\ga)\leq -\ep \ \ \IF y_\ga\in\gO,\\
\label{st-visco3}
\min\{H_i(p_\ga),\,F(X_\ga, p_\ga, \tilde v(x_\ga),x_\ga)\}\leq -\ep\rho \ \ \IF x_\ga\in\pl_i\gO,\,i=1,2, \\
\label{st-visco4}
\min\{H_i(p_\ga),\,F(Y_\ga, p_\ga, \tilde w(y_\ga),y_\ga)\}\leq -\ep\rho \ \ \IF y_\ga\in\pl_i\gO,\,i=1,2. 
\end{gather}

It is a standard observation (see, e.g.,  \cite{CIL1992}) that for some 
$x_0\in\ol\gO$, as $\ga \to \infty$ along a sequence $\{\ga_j\}$,
\[
(x_\ga,y_\ga) \to (x_0,x_0),  \quad 
\tilde v(x_\ga) \to \tilde v(x_0), \quad
\tilde w(y_\ga)\to \tilde w(x_0), \quad \AND \quad 
\ga^{-1}f(\ga(x_\ga-y_\ga)) \to 0. 
\]  
In particular, \erf{f2} ensures that, as $\ga=\ga_j$ and $j\to\infty$,
\beq\label{est-(x-y)}
\ga|x_{\ga}-y_\ga|^2 \ \to \ 0.
\eeq
Furthermore, plugging $x=\xi+q|\xi|$, with $|q|=1$, and $\xi=\ga(x_\ga-y_\ga)$ into \erf{semi-cv} yields
\[
\ga|x_\ga-y_\ga||p_\ga|\leq f(\ga(x_\ga-y_\ga))+B\ga^2|x_\ga-y_\ga|^2,
\]  
which, combined with \erf{est-(x-y)} and the previous observations, shows that, as $\ga=\ga_j$ and $j\to\infty$, 
\beq \label{est-p}
\ga^{-1}|p_\ga||x_\ga-y_\ga| \,\to\, 0.
\eeq

Focusing on $\ga=\ga_j$, with large $j$, we may assume that 
$\tilde v(x_\ga)>\tilde w(y_\ga)$.  
Observe by the choice of $f$ that if $x_\ga\in\pl_1\gO$, then 
$(x_\ga-y_\ga)\cdot e_1=-y_\ga\cdot e_1\leq 0$ and, by \erf{f1}, 
\[H_1(p_\ga)=H_1(Df(\ga(x_\ga-y_\ga))\geq 0.
\]
Similarly, if $x_\ga\in\pl_2\gO$, then $(x_\ga-y_\ga)\cdot e_2= -y_\ga\cdot e_2\leq 0$ and $H_2(p_\ga)=H_2(Df(x_\ga-y_\ga))\geq 0$.
Also, if $y_\ga\in\pl_1\gO$, then 
$(x_\ga-y_\ga)\cdot e_1=\ga x_\ga\cdot e_1\geq 0$ and $H_1(p_\ga)
=H_1(Df(x_\ga-y_\ga))\leq 0$. 
If $y_\ga\in\pl_2\gO$, then $(x_\ga-y_\ga) \cdot e_2 =x_\ga\cdot e_2\geq 0$ 
and $H_2(p_\ga)=H_2(Df(x_\ga-y_\ga))\leq 0$. 
Taking into account these sign properties, we obtain by \erf{st-visco1}--\erf{st-visco4}, 
\[\bald
 &F(X_\ga,p_\ga,\tilde v(x_\ga),x_\ga)\leq -\ep\rho,  
\\&F(Y_\ga,p_\ga,\tilde w(y_\ga),y_\ga)\geq \ep\rho,  
\eald
\]
Since $u\mapsto F(X,p,u,x)$ is increasing on $\R$ by \erf{F1}, we have 
\[
F(X_\ga,p_\ga,\tilde w(y_\ga),x_\ga)\leq -\ep\rho,
\]
and
\[
-2\ep\rho\geq F(X_\ga,p_\ga,\tilde w(y_\ga),x_\ga)-F(Y_\ga,p_\ga,\tilde w(y_\ga),y_\ga),
\]
which yields, together with \erf{F3},
\[
-2\ep\rho\geq -\go(3B\ga|x_\ga-y_\ga|^2+|x_\ga-y_\ga|(|p_\ga|+1)).
\]
Recalling \erf{est-(x-y)} and \erf{est-p} and sending $j\to\infty$ for $\ga=\ga_j$, we 
get a contradiction, $\go(0)=0\leq -2\ep\rho$. The proof is complete.  
\eproof

\bye